%% file: author.tex
\documentclass{article}

\usepackage{graphicx}
\usepackage{ucs}
\usepackage[utf8x]{inputenc}
\usepackage{amsmath}
\usepackage{amsfonts}
\usepackage[english]{babel}
\usepackage{fontenc}
\usepackage{lmodern}

\DeclareMathOperator{\spann}{span}
\DeclareMathOperator{\divergence}{div}
\newcommand{\R}{\mathbb{R}}
\usepackage{tikz}
\usetikzlibrary{calc}

\newtheorem{theorem}{Theorem}

\begin{document}

\title{\Large The Influence of Quadrature Errors \\on Isogeometric Mortar Methods}

\author{
{\sc Ericka Brivadis\thanks{Corresponding author. Email: ericka.brivadis@iusspavia.it}
 }\\[2pt]
 Istituto Universitario di Studi Superiori Pavia,\\ Palazzo del Broletto, Piazza della Vittoria 15, 27100 Pavia, Italy\\[6pt]
 {\sc and}\\[6pt]
  {\sc Annalisa Buffa\thanks{Email: annalisa@imati.cnr.it} }\\[2pt]
  Istituto di Matematica Applicata e Tecnologie Informatiche del CNR,\\ Via Ferrata 1, 27100 Pavia, Italy\\[6pt]
 {\sc and}\\[6pt]
{\sc Barbara Wohlmuth\thanks{Email: wohlmuth@ma.tum.de}  and
Linus Wunderlich\thanks{Email: linus.wunderlich@ma.tum.de}
}\\[2pt]  M2 -- Zentrum Mathematik, Technische Universit\"at M\"unchen,\\ Boltzmannstra\ss{}e 3, 85748 Garching, Germany}
\date{Oktober 1, 2014}
\maketitle
\begin{abstract}
Mortar methods have recently been shown to be well suited for isogeometric analysis. We review the recent mathematical analysis and then investigate the variational crime introduced by quadrature formulas for the coupling integrals. Motivated by finite element observations, we consider a quadrature rule purely based on the slave mesh as well as a method using quadrature rules based on the slave mesh and on the master mesh, resulting in a non-symmetric saddle point problem. While in the first case reduced convergence rates can be observed, in the second case the influence of the variational crime is less significant.
\end{abstract}

\section{Introduction}
\label{sec:1}
\input{01_introduction.tex}

\section{Isogeometric Mortar Methods}
\label{sec:2}
\input{02_iga_mortar.tex}

\section{Mortar Integrals}
\label{sec:3}
\input{03_quadrature_review}

\section{Numerical Results}
\label{sec:4}
\input{04_numerics.tex}

\subsubsection{Slave Integration Approach}\label{sec:slave_int}
\input{04_1_numerics.tex}
\subsubsection{Non-symmetric Approach}
\input{04_2_numerics.tex}

\subsection{Three-dimensional  Example}
\input{04_3_numerics.tex}

\section{Conclusion}
\input{05_conclusion.tex}

\section*{acknowledgement}
The first author has been supported by Michelin under the contract A10-4087.
The third and the fourth authors have been supported by the International Research Training Group IGDK 1754, funded by the German Research Foundation (DFG) and the Austrian Science Fund (FWF), and by the German Research Foundation (DFG) in Project WO 671/15-1 and the Priority Programme “Reliable Simulation Techniques in Solid Mechanics. Development of Non-standard Discretisation Methods, Mechanical and Mathematical Analysis” (SPP 1748). 
The last author was supported by the Elite Network of Bavaria through its graduate program TopMath and the TUM Graduate School through its Thematic Graduate Center TopMath.
All supports are gratefully acknowledged.

 \bibliographystyle{alpha}
 \bibliography{bibliography}

 \end{document}

%% file: 01_introduction.tex
Isogeometric analysis, introduced in 2005 by  Hughes et al. in~\cite{hughes:05}, is a family of  methods that use B-splines and non-uniform rational B-splines (NURBS) as basis functions to construct numerical approximations of partial differential equations (PDEs). 
With isogeometric methods, the computational domain is generally split into patches. Within this framework, techniques to couple the numerical solution on different patches are required. To retain the flexibility of the meshes at the interfaces, mortar methods are very attractive.

Mortar  methods are a popular tool for the coupling of non-matching meshes, originally introduced for spectral and finite element methods~\cite{bernardi:94,ben_belgacem:97,ben_belgacem:99}. They were successfully applied in the context of isogeometric analysis~\cite{hesch:12,bletzinger:14,dornisch:14}.
A mathematical analysis enlightening the use of different dual spaces was recently presented in~\cite{brivadis:14}. 
In this paper, starting from these latter results, we focus on one particular challenge in the realization of a mortar method, namely, the evaluation of the interface integrals which contain a product of functions defined on non-matching meshes. 

This article is structured as follows. In Sec.~\ref{sec:2}, we recall the basics of isogeometric mortar methods. In Sec.~\ref{sec:3}, we consider a review of numerical quadrature for mortar integrals such as additional aspects in the case of isogeometric analysis, illustrated by numerical results in Sec.~\ref{sec:4}.

%% file: 02_iga_mortar.tex
In the following, we briefly present isogeometric mortar methods, for more details we refer to~\cite{brivadis:14}. 
After stating the problem setting, we review isogeometric parametrizations, describe the domain decomposition into several NURBS patches and finally discuss suitable coupling spaces.

Let $\Omega \subset \R^d$, $d$ the dimension being $2$ or $3$, be a bounded domain, $\alpha,\, \beta \in L^{\infty}(\Omega)$, $\alpha > \alpha_0 > 0$ and $\beta \geq 0$. We consider the following  second order elliptic boundary value problem with homogeneous Dirichlet conditions
\begin{subequations}\label{eq:strong_formulation}
\begin{align} 
		-\divergence (\alpha\nabla u) + \beta u &= f \quad \text{ in } \Omega,\\
		u &= 0 \quad \text{ on } \partial\Omega_D = \partial\Omega.
\end{align}
\end{subequations}
We assume $\alpha$ and $\beta$ to be piecewise sufficiently smooth.

\subsection{Isogeometric Parametrization} \label{sec:part_iga}
Here, we present  isogeometric concepts and notations used throughout the paper, and refer to the classical literature~\cite{piegl:97,bazilevs:06,Schumaker:07,hughes:09} for more details.

 Let us denote by $p$ the degree of the univariate B-splines and by $\Xi $ an open uni\-variate knot vector, where the first and last entries are repeated $(p+1)$-times, i.e.,
\[
	\Xi = \{ 0 = \xi_1 = \ldots = \xi_{p+1} < \xi_{p+2} \leq \ldots \leq \xi_{n} < \xi_{n+1} = \ldots= \xi_{n+p+1} = 1\}.
\]
Let us define $Z = \{\zeta_1,\, \zeta_2,\, \ldots,\, \zeta_{E}\}$ as the knot vector without any repetition, also called breakpoint vector.
For each breakpoint $\zeta_j$ of $Z$, we define its multiplicity $m_j$ as its number of repetitions in $\Xi$.  The Cox-de~Boor algorithm, see~\cite{hughes:09}, defines $n$ univariate B-splines $\widehat{B}_i^p(\zeta)$, $i=1,\ldots, n$, based on the univariate knot vector $\Xi$ and the degree $p$. We denote by $S^p(\Xi)=\spann\{\widehat{B}_i^p(\zeta), \,i=1,\, \ldots,\, n\}$ the corresponding spline space. The smoothness of B-splines is defined by the breakpoint multiplicity. More precisely, each basis function is $C^{p-m_j}$ at each $\zeta_j \in Z$.

To define multivariate spline spaces, we introduce the multivariate knot vector $ \mathbf{\Xi} =(\Xi_1 \times \Xi_2 \times \ldots \times \Xi_d)$, and  for simplicity of notations assume in the following that the degree is the same in all parametric directions and denote it by $p$. Multivariate B-splines $\widehat{B}_{\mathbf{i}}^p({\boldsymbol \zeta})$ are defined by tensor product of univariate B-splines for each multi-index $\mathbf i \in \mathbf{I} =\{(i_1, \, \ldots, \, i_d): 1\leq i_\delta \leq n_{\delta}\}$. We denote by $S^p(\mathbf{\Xi})$ the corresponding spline space in the parametric domain. 

Given a set of positive weights $\omega_{\mathbf{i}}$, we define NURBS functions $\widehat{N}_{\mathbf{i}}^p({\boldsymbol \zeta })$ as rational functions of B-splines and the weight function $\widehat{W} = \sum_{ \mathbf{i} \in \mathbf{I}} \omega_{\mathbf{i}} \,\widehat{B}_{\mathbf{i}}^p({\boldsymbol \zeta})$, and set $N^p(\mathbf{\Xi})$ as the multivariate NURBS space in the parametric domain.

For a set of control points $\mathbf{C}_{\mathbf{i}} \in \mathbb{R}^d$, $\mathbf{i} \in \mathbf{I}$, we define a parametrization of a NURBS surface ($d=2$) or solid ($d=3$) as a linear combination of NURBS and  control points
\begin{equation*}
	\mathbf{F}({\boldsymbol \zeta})=\displaystyle \sum_{\mathbf{i} \in \mathbf{I}} \mathbf{C}_{\mathbf{i}} \,\widehat{N}_{\mathbf{i}}^p({\boldsymbol \zeta }),
\end{equation*}
and assume the regularity stated in~\cite[Assumption 3.1]{beirao:14}.

The knot vector ${\mathbf \Xi}$ forms a mesh in the parametric space $\widehat \Omega$. We define the physical mesh $\mathcal{M}$ as the image  of this parametric mesh through $\mathbf{F}$, and denote by $\mathbf{O}$ its elements.
The $h$-refinement procedure, see~\cite[Section 2.1.3]{beirao:14}, yields a family of meshes denoted $\mathcal{M}_{h}$, each mesh being a refinement of the initial one. We assume quasi-uniformity for each mesh.  

\subsection{Description of the Computational Domain}

Let the domain $\Omega$ be decomposed into $K$ non-overlapping domains $\Omega_k$, i.e.,
\[
\overline{\Omega} = \bigcup_{k=1}^K \overline{\Omega}_k, \text{ and } \Omega_i \cap \Omega_j = \emptyset, i \neq j.
\]

Each subdomain is a NURBS geometry, i.e., there exists a NURBS parametrization $\mathbf F_k$ based on a knot vector $\mathbf{\Xi}_k$ and a degree $p_k$, see Sec.~\ref{sec:part_iga}, such that $\Omega_k$ is the image of the parametric space $\widehat \Omega = (0,1)^d$ by $\mathbf F_k$.

For $1\leq k_1,\, k_2 \leq K $, $k_1\neq k_2$, we define the interface as the interior of the intersection of the boundaries, i.e., $\overline{\gamma}_{k_1k_2} = \partial {\Omega}_{k_1} \cap \partial {\Omega}_{k_2}$, where ${\gamma}_{k_1k_2}$ is open. Let the non-empty interfaces be enumerated by $\gamma_l$, $l = 1,\,\ldots,\, L$,  and let us define the skeleton $\Gamma = \bigcup_{l=1}^L \gamma_l$ as the union of all interfaces. For each interface, one of the adjacent subdomains is chosen as the master side and one as the slave side. This choice is arbitrary but fixed. 
We denote the index of the former by $m(l)$, the index of the latter one by ${s(l)}$, and thus $\overline{\gamma}_l=\partial \Omega_{m(l)}\cap \partial \Omega_{s(l)}$. On the interface $\gamma_l$, we define the outward normal ${\bf  n}_l$ of the master side $\partial \Omega_{m(l)}$ and denote  by $\displaystyle{\frac{\partial u}{{\partial\bf n}_l}}$ the normal derivative on $\gamma_l$ from the master side.

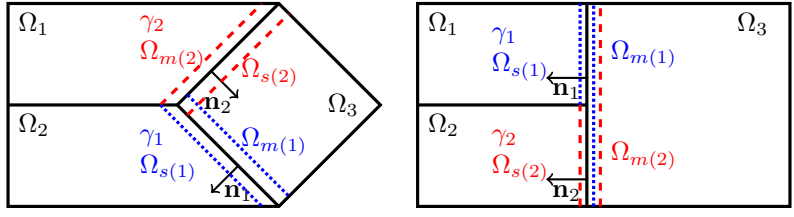
\begin{figure}
	\begin{center}
	\input{figures/mortar_setting_picture}
	\end{center}
	\caption{Geometrical conforming case (left) and  slave conforming case (right).}
	\label{mortar:mortar_setting} 
\end{figure}

We assume that for each interface the pull-back with respect to the slave domain is a whole face of the unit $d$-cube in the parametric space, which we call a \emph{slave conforming} situation, see the right setting in Fig.~\ref{mortar:mortar_setting}.
If we also assume that the pull-back with respect to the master domain is a whole face of the unit $d$-cube, we are in a fully geometrically conforming situation, see the left picture of Fig.~\ref{mortar:mortar_setting}.

For each $\Omega_k$, we introduce the space $H^1_*(\Omega_k)=\{ v_k \in H^1(\Omega_k),  v_{k |_{\partial \Omega \cap \partial \Omega_k}}=0\}$, where we use standard Sobolev spaces, as defined in~\cite{grisvard:11}, endowed with their usual norms. 
For any interface $\gamma_l \subset \partial \Omega_{s(l)}$, we define by $H^{1/2}_{00}(\gamma_l) \subset H^{1/2}(\partial \Omega_{s(l)})$ the space of all functions that can be trivially extended on $\partial \Omega_{s(l)} \setminus \gamma_l$ by zero to an element of $H^{1/2}(\partial \Omega_{s(l)})$. Note that $H^{1/2}(\partial \Omega_{s(l)})$ is the trace space of $H^1(\Omega_{s(l)})$. The dual space of $H^{1/2}_{00}(\gamma_l)$ is denoted $H^{-1/2}(\gamma_l)$.
In order to set a global functional framework on $\Omega$, we consider the broken Sobolev spaces $V= \prod_{k=1}^K H^1_*(\Omega_k)$, endowed with the broken norm $ \| v \|_{V}^2 = \sum_{k=1}^K \| v \|_{H^1(\Omega_k)}^2$, and $M= \prod_{l=1}^L H^{-1/2}(\gamma_l)$.

The mortar method is based on a weak coupling between different subdomains. Each subdomain is discretized independently and a weak coupling is performed on each interface. 
From now on, we assume that jumps of $\alpha$ and $\beta$ are solely located at the skeleton, and we define the  linear and bilinear forms $a\colon V \times V \rightarrow \mathbb{R}$ and  
 $f\colon V \rightarrow \mathbb{R}$, such that  
 \[a(u,v) = \sum_{k=1}^K \int_{\Omega_k} 
\alpha \nabla u \cdot  \nabla v +  \beta\, u\, v ~ \mathrm{d}{\bf x}, \quad
f(v) =\sum_{k=1}^K \int_{\Omega_k} f v ~ \mathrm{d}{\bf x}.
\]
We remark that the standard weak formulation of~\eqref{eq:strong_formulation}, where no weak coupling is necessary and which is uniquely solvable, reads as follows: Find $u\in H^1_0(\Omega)$, such that
\begin{align} \label{eq:weak_formulation}
a(u,v) = f(v), \quad v\in H_0^1(\Omega).
\end{align}

\subsection{Isogeometric Mortar Discretization}

In the following, we set our non-conforming approximation framework. On each subdomain $\Omega_k$, based on the NURBS parametrization, we introduce the approximation space $V_{k,h}=\{v_k=\widehat{v}_k \circ \mathbf{F}_k^{-1}, \widehat{v}_k \in N^{p_k}(\mathbf{\Xi}_k) \}$. We recall that under the assumptions on the mesh $\mathcal{M}_{k,h}$ and on the parametrization $\mathbf{F}_k$, this NURBS space has optimal approximation properties, see, e.g.,~\cite{bazilevs:06}. On $\Omega$, we define the discrete product space $V_h = \prod_{k=1}^K V_{k,h} \subset V$, which forms a $H^1(\Omega)$ non-conforming space discontinuous over the interfaces. We denote in the following the maximal mesh size $h = \max_k h_k$ as the mesh parameter.

On the skeleton $\Gamma$, we define the discrete product Lagrange multiplier space $M_h$ as $M_h = \prod_{l=1}^L M_{l,h}$, where $M_{l,h}$ denotes one of the two following choices which were shown to be well suited in~\cite{brivadis:14}.
The first choice $M_{l,h}^0$ is the spline space of degree $p_{s(l)}$, defined on the interface $\gamma_l$ based on the interface knot vector of the slave body $\Omega_{s(l)}$. Note that in the presence of any cross point, a suitable modification, e.g., a local degree reduction as presented in~\cite[Section 4.3]{brivadis:14}, has to be applied. 
The alternative choice $M_{l,h}^{2}$ is an order $(p_{s(l)}-2)$ spline space defined on the interface $\gamma_l$ based on the interface knot vector of the slave body $\Omega_{s(l)}$ for which the definition requires the trace space of $V_{s(l),h}$ to be a subset of $C^1(\gamma_l)$. 
Moreover, a third dual space $M_{h,l}^1$ which is a spline space of degree $(p_{s(l)}-1)$, could be considered. This choice is not inf-sup stable, so we do not consider it any further.

The saddle point formulation of the isogeometric mortar method reads as follows: 
Find $ (u_h, \lambda_h) \in V_h \times M_h,$  such that 
\begin{subequations}\label{eq:discrete_spp}
\begin{align}
		a(u_h, v_h)+ b(v_h, \lambda_h) &= f(v_h), \quad   v_h \in V_h,\\
		b(u_h, \mu_h) &= 0, \quad  \mu_h \in M_h,
\end{align}
\end{subequations}
where $b(v,\mu) = \sum_{l=1}^L \int_{\gamma_l} \mu [v]_l ~\mathrm{d}\sigma$ and  $[\cdot]_l$ denotes the jump from the master to the slave side over $\gamma_l$.
We note that the Lagrange multiplier $\lambda_h$ gives an approximation of the normal flux across the skeleton.

It is well known from the theory of mixed and mortar methods that two abstract requirements on each interface guarantee the method to be well-posed and  of optimal order, see~\cite{ben_belgacem:99}. Namely, an appropriate approximation order of the dual space and a uniform inf-sup stability between the primal space and the dual space. 
Note that for simplicity of notations, we assume the same type of dual space to be used for all interfaces.
The following theorem is shown in~\cite{brivadis:14} and guarantees a-priori bounds.
\begin{theorem} \label{thm:convergence_rates}
Let $\theta = 0$ if $M_{h,l} = M_{h,l}^0$ and $\theta = 1/2$ if   $M_{h,l} = M_{h,l}^2$.
For $u\in H^{\sigma+1}(\Omega)$, $1/2<\sigma \leq \min_{k,l} (p_k-\theta)$, the solution of~\eqref{eq:weak_formulation} and $(u_h, \lambda_h)$ the non-conforming approximation, see~\eqref{eq:discrete_spp},  it holds
\[
	\frac{1}{h^2} \|u-u_h\|_{L^2(\Omega)}^2 + \| u - u_h \|_{V}^2  \leq C \sum_{k=1}^K h^{2\sigma}_k \|u \|^2_{H^{\sigma+1}(\Omega_k)}
\]
for the primal solution and
\[
\sum_{l=1}^L \| \alpha \frac{\partial u}{\partial \mathbf{n}_l} - \lambda_h \|_{H^{-1/2}(\gamma_l)}^2   \leq C \sum_{k=1}^K h^{2\sigma}_k \|u \|_{H^{\sigma+1}(\Omega_k)}^2
\]
for the dual solution.
With $0 < C < \infty$ a generic constant that is independent of the mesh sizes but possibly depends on $p_k$.
\end{theorem}
We highlight that while for both pairings the inf-sup stability is satisfied, the approximation order of the lower order dual space $M_{l,h}^2$ is only close to optimal. Indeed, in this case an order of $\sqrt{h}$ is lost in the convergence order.  

%% file: figures/mortar_setting_picture.tex
\begin{tikzpicture}[scale=0.9]

\tikzstyle{style1} = [black, very thick]
\tikzstyle{style2} = [red, very thick, dashed]
\tikzstyle{style3} = [blue, very thick, dash pattern=on \pgflinewidth off 1pt]

\def\a{5.5};
\def\b{2.5};
\def\c{4};
\def\d{2.5};
\def\h{3};
\def\hh{1.5};
\def\f{0.1};
\def\s{0.1};
\def\ss{0.2};

\coordinate (0) at (0, 0);
\coordinate (1) at (\c, 0);
\coordinate (2) at (\a, \hh);
\coordinate (3) at (\c, \h);
\coordinate (4) at (0, \h);
\coordinate (5) at (0, \hh);
\coordinate (6) at (\b, \hh);

\coordinate (7) at ($(6)!0.1!(5)$);
\coordinate (8) at ($(6)!0.1!(3)$);
\coordinate (9) at ($(6)!0.1!(1)$);

\coordinate (10) at ($(7) + (3) - (6)$);
\coordinate (11) at ($(7) + (1) - (6)$);
\coordinate (12) at ($(8) + (1) - (6)$);
\coordinate (13) at ($(9) + (3) - (6)$);

\def\ln{0.4};
\coordinate(14)at (\c-0.6,0.6);
\coordinate(15)at(\c-0.6-\ln,0.6-\ln);
\draw[thick,->]  (14) -- (15);
\coordinate(18)at(\c-\ln-0.2,0.2);
\node[style1] at ($(18)$) {${\bf n}_1$};

\coordinate(16)at (\b+0.5,\hh+0.5);
\coordinate(17)at(\b+0.5+\ln,\hh+0.5-\ln);
\draw[thick,->]  (16) -- (17);
\coordinate(19)at(\b+0.5+0.1,\hh+0.5-\ln-0.1);
\node[style1] at ($(19)$) {${\bf n}_2$};

\draw[style1] (0) -- (1) -- (2) -- (3) -- (4) -- cycle;
\draw[style1]  (5) -- (6) -- (3);
\draw[style1]  (6) -- (1);

\draw[style2] (7) -- (10);
\draw[style3] (7) -- (11);
\draw[style3] (8) -- (12);
\draw[style2] (9) -- (13);

\node[right, style1] at ($(0)!0.9!($(0)+(0,\h)$)$) {$\Omega_1$};
\node[right, style1] at ($(0)!0.4!($(0)+(0,\h)$)$) {$\Omega_2$};
\node[style1] at ($(5)!0.9!(2)$) {$\Omega_3$};

\node[style2, text width=1cm] at ($(6) + (0, 0.31*\h)$) {$\gamma_2$\\$\Omega_{m(2)}$};
\node[style3, text width=1cm]  at ($(6) - (0, 0.25*\h)$) {$\gamma_1$\\$\Omega_{s(1)}$};

\node[style3, text width=1cm] at ($(1)+(0, 0.32*\h)$) {$\Omega_{m(1)}$};
\node[style2, text width=1cm] at ($(1)+(0, 0.66*\h)$) {$\Omega_{s(2)}$};

\begin{scope}[xshift=1.1*\a cm]
\coordinate (0) at (0, 0);
\coordinate (1) at (\a, \h);
\coordinate (2) at (\d, 0);
\coordinate (3) at (\d, \h);
\coordinate (4) at (0, \hh);
\coordinate (5) at (\d, \hh);
\coordinate (6) at ($(2) - (\s, 0)$);
\coordinate (7) at ($(2) + (\s, 0)$);
\coordinate (8) at ($(2) + (\ss, 0)$);
\coordinate (9) at ($(5) - (\s, 0)$);
\coordinate (10) at ($(5) - (\s, 0)$);
\coordinate (11) at ($(3) - (\s, 0)$);
\coordinate (12) at ($(3) + (\s, 0)$);
\coordinate (13) at ($(3) + (\ss, 0)$);
\coordinate (14) at ($(0)!0.3!($(0)+(\a, 0)$)$);
\coordinate (15) at ($(0)!0.62!($(0)+(\a, 0)$)$);

\def\lnb{0.6};
\coordinate(16)at (\d,\hh+0.4);
\coordinate(17)at(\d-\lnb,\hh+0.4);
\coordinate(18)at(\d-\lnb/2,\hh+0.2);
\draw[thick,->]  (16) -- (17);
\node[style1] at ($(18)$) {${\bf n}_1$};

\coordinate(19)at (\d,0.4);
\coordinate(20)at(\d-\lnb,0.4);
\coordinate(21)at(\d-\lnb/2,0.2);
\draw[thick,->]  (19) -- (20);
\node[style1] at ($(21)$) {${\bf n}_2$};

\draw[style1]  (0) rectangle (1);
\draw[style1]  (2) -- (3);
\draw[style1]  (4) -- (5);
\draw[style2] (6) -- (9);
\draw[style3] (10) -- (11);
\draw[style3] (7) -- (12);
\draw[style2] (8) -- (13);

\node[right, style1] at ($(0)!0.9!($(0)+(0,\h)$)$) {$\Omega_1$};
\node[right, style1] at ($(0)!0.4!($(0)+(0,\h)$)$) {$\Omega_2$};
\node[style1] at ($(0)!0.9!(1)$) {$\Omega_3$};

\node[style3, text width=1cm, ] at ($(14)+(0, 0.75*\h)$) {$\gamma_1$\\$\Omega_{s(1)}$};
\node[style2, text width=1cm, ] at ($(14)+(0, 0.25*\h)$) {$\gamma_2$\\$\Omega_{s(2)}$};

\node[style3, text width=1cm, ] at ($(15)+(0, 0.75*\h)$) {$\Omega_{m(1)}$};
\node[style2, text width=1cm, ] at ($(15)+(0, 0.25*\h)$) {$\Omega_{m(2)}$};

\end{scope}

\end{tikzpicture}

%% file: 03_quadrature_review.tex
To evaluate the bilinear form $b(v,\mu)$, we need to evaluate for each interface $\gamma_l$ the mortar integrals $\int_{\gamma_l} \mu\,  v^+  \mathrm{d}\sigma$ and $\int_{\gamma_l} \mu\,  v^-  \mathrm{d}\sigma$, where $v^+$ denotes the trace of $v$ from the master domain $\Omega_{m(l)}$ and $v^-$ the trace of $v$ from the slave domain $\Omega_{s(l)}$. To simplify the notation, let us restrict ourselves to the case of one single interface and drop the index $l$ in the following.

One particular challenge in the realization of a mortar method is the evaluation of the first interface integral due to the product $\mu\, v^+$ of functions which are defined on non-matching meshes, see \cite{bertoluzza:04} for a method to bypass it in a finite element/wavelet context.  
Any quadrature rule based on the slave mesh does not respect the mesh lines of the master mesh and vice versa for a quadrature based on the master mesh. 

It is obvious that the use of a suitable quadrature rule based on a merged mesh, i.e., a mesh which respects the reduced smoothness of the master and slave functions at their respective mesh lines, leads to an exact evaluation of the integral.
However, the construction of this auxiliary mesh commonly named segmentation process is challenging, especially in the three dimensional case since the shape of the elements is not unique and difficult to determine, see, e.g.,~\cite{mcdevitt:00,pusoa:04,pusob:04,hesch:12,dittmann:14}. 
 Note that  in an isogeometric context the merged mesh needs to be constructed in the physical space and then pulled back to the parametric space for each subdomain. 
The complexity of constructing such a mesh becomes even more severe in the case of non-linear and time-dependent problems, where the relative position of the meshes changes in every time or load step which implies to recompute the merged mesh at every step. 

Due to this computational complexity, it has been seen very appealing to use a higher order quadrature rule either based on the slave mesh or on the master mesh, see~\cite{fischer:05,tur:09,delorenzis:11} for some applications in finite element and isogeometric analysis contexts. However in the finite element case, early results in~\cite{maday:97, wohlmuth:02} showed that this strategy does not yield optimal methods. More precisely, in the case the master mesh is chosen, the best approximation error is affected, while in contrast in the case the slave mesh is chosen it is the consistency error. Numerical results confirmed the lack of optimality with the master integration approach, while with the slave integration approach reasonable results were obtained although not optimal in terms of the Lagrange multiplier norm.

Due to the global smoothness of splines, one could expect the sensitivity with respect to the quadrature rules for isogeometric methods to be less than for finite element methods. 
In the mortar context, according to the finite element results, it seems interesting to consider a slave integration rule. And, in case of maximal regularity, i.e., $V_{k,h}\subset C^{p_k-1}(\Omega_k)$ one also might expect the
quadrature error on a non-matching mesh to be significantly smaller than in the finite element case. 

Let us denote the quadrature rule based on the boundary mesh of the slave domain as $\sum_-$, i.e.,
$
\int_\gamma \mu v^+ \mathrm{d}\sigma \approx {\sum}_- \mu v^+.
$
We precise that in the examples a Gaussian quadrature rule with a various number of points is used. 
The mortar method with pure slave integration is obtained by evaluating all interface integrals in~\eqref{eq:discrete_spp} using this quadrature rule, i.e., the discrete system reads as follows: Find $(\widetilde{u}_h, \widetilde{\lambda}_h)\in V_h\times M_h$, such that
\begin{align*}
		a(\widetilde{u}_h, v_h)+  \sum\nolimits_- (v_h^+ - v_h^-) \widetilde{\lambda}_h &= f(v_h), \quad   v_h \in V_h,\\
		 \sum\nolimits_-  (\widetilde{u}_h^+ - \widetilde{u}_h^-) \mu_h &= 0, \quad  \mu_h \in M_h.
\end{align*}
The notation $\widetilde{\cdot}$ is used to stress the difference to the discrete solution with exact integration.

In the next section, we present numerical examples which show severe disturbances even in the isogeometric case. Hence, even though the global smoothness of the integrated function is increased compared to the finite element case, a non-matching integration approach reduces the convergence order drastically.

Moreover, we consider an alternative approach which was proposed in~\cite{maday:97, wohlmuth:02} using both integration rules. Additionally denoting $\sum_+$ a quadrature rule based on the boundary mesh of the master domain $\Omega_m$, this approach which results in a non-symmetric saddle point problem, reads as follows: Find $(\widetilde{u}_h, \widetilde{\lambda}_h)\in V_h\times M_h$, such that
\begin{align*}
		a(\widetilde{u}_h, v_h)+  \sum\nolimits_+ v_h^+ \widetilde{\lambda}_h - \sum\nolimits_- v_h^- \widetilde{\lambda}_h &= f(v_h), \quad   v_h \in V_h,\\
		 \sum\nolimits_-  (\widetilde{u}_h^+ - \widetilde{u}_h^-) \mu_h &= 0, \quad  \mu_h \in M_h.
\end{align*}
The non-symmetric saddle point problem{, which corresponds to a Petrov--Galerkin approach in the primal formulation,} was motivated by different requirements for the integration of the primal and dual test functions. Numerical examples showed error values very close to the case of exact integration, but we note that from the theoretical side even the well-posedness of the non-symmetric saddle point problem remains unclear. In the next section, we present numerical examples which show that also in an isogeometric context the results are generally close to those from the exact integration case.

%% file: 04_numerics.tex
In this section, we consider two-dimensional and three-dimensional settings in order to observe the effects of inexact quadrature rules, as presented in Sec.~\ref{sec:3}, on the optimality of the mortar method. We first set the problem settings, and then give the results of several studies.  

\subsection{Two-dimensional Example}
\begin{figure}
\begin{center}
\includegraphics[height=0.41\textwidth]{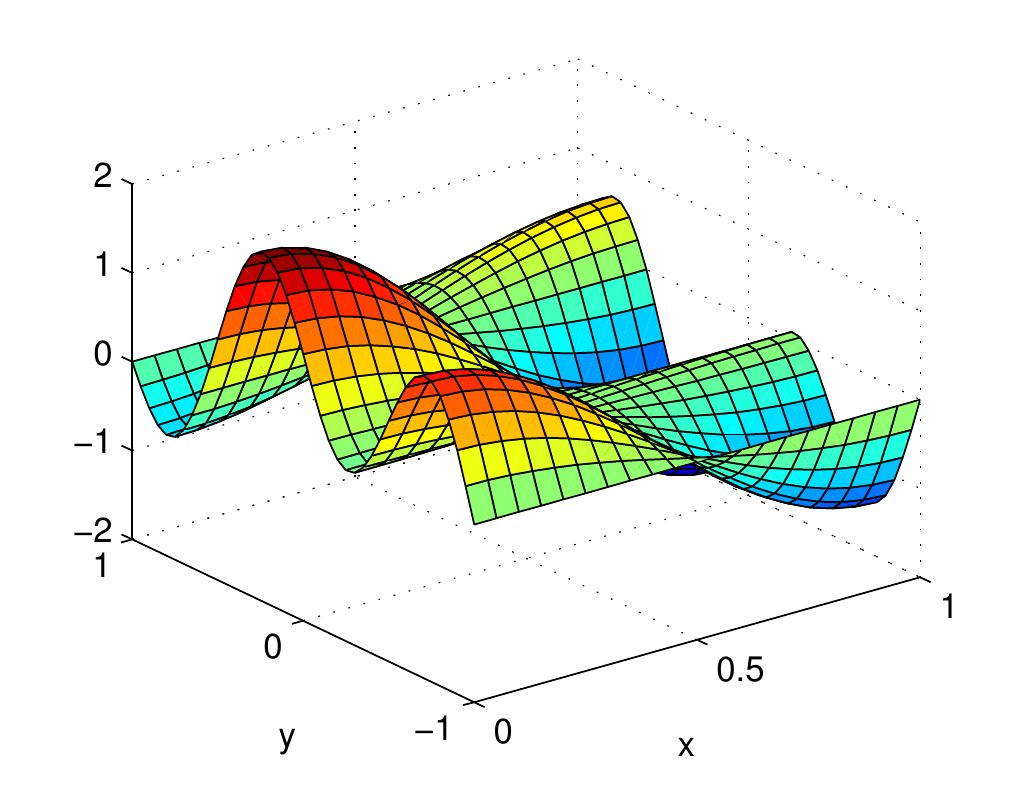}
\includegraphics[height=.265\textwidth]{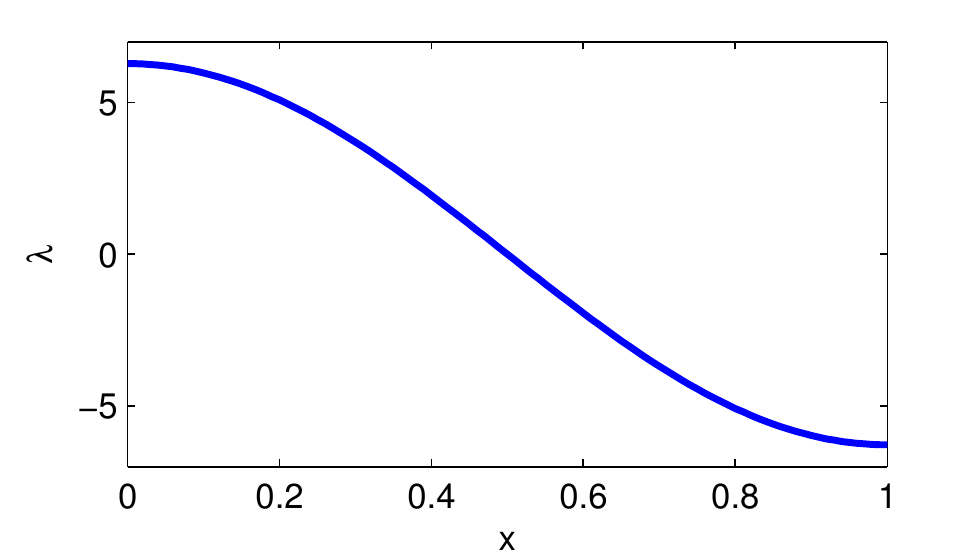}
\caption{Left: Primal solution on $\Omega$. Right: Lagrange multiplier along the interface.}\label{fig:solutions}
\end{center}
\end{figure}

As a first example, let us consider the Poisson problem $-\Delta u=f$ solved on the domain $\Omega=(0,1)\times (-1,1)$ which is decomposed into two patches by the interface $\gamma = \{(x,y) \in \Omega,\, y=0\}$. The upper domain is set as the slave domain. The internal load and the boundary conditions are manufactured to have the analytical solution 
\[
u(x,y) = \cos\left( \pi x \right)  (\cos\left( \frac{\pi}{2} y \right)+\sin\left( 2\pi y \right)).
\] 
The normal derivative on the interface is given by  { ${\partial u}/{\partial {\bf n} }(x) = 2 \pi \cos\left( \pi x \right),$ }
 see  Fig.~\ref{fig:solutions}. Neumann conditions are applied on the left and right boundary parts, such that no cross point modification is necessary.

\begin{figure}
\begin{center}
\includegraphics[height=.43\textwidth]{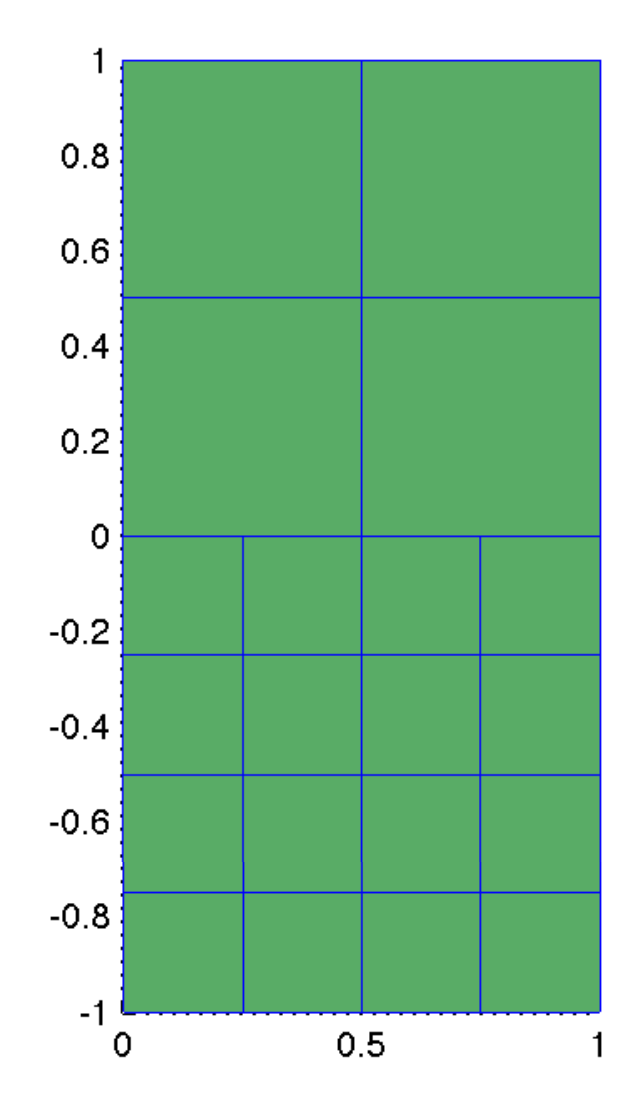}
\includegraphics[height=.43\textwidth]{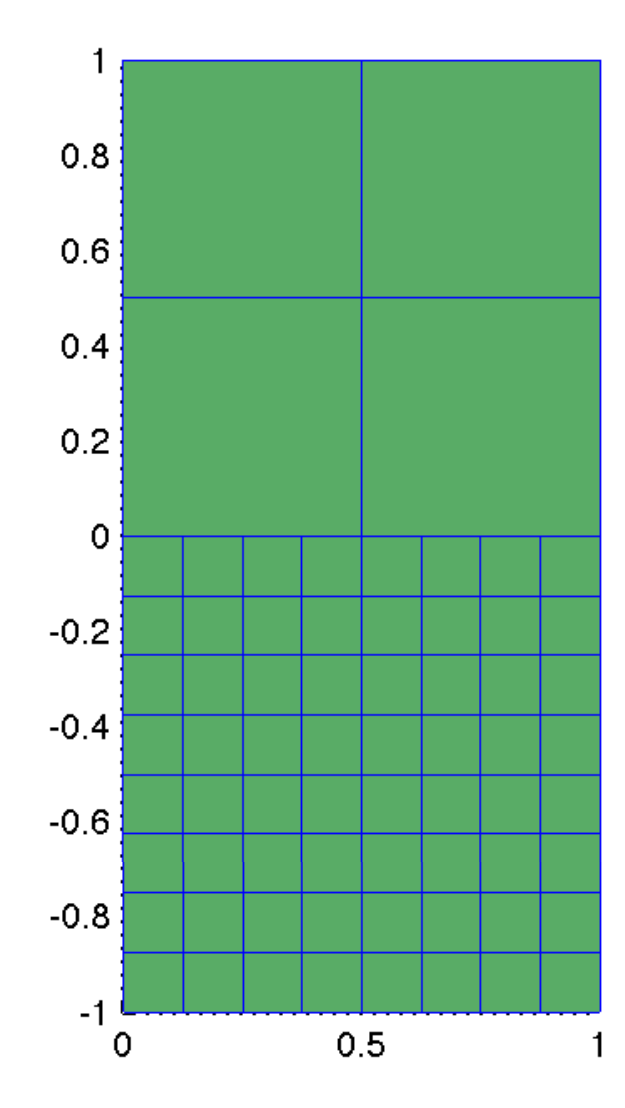}
\includegraphics[height=.43\textwidth]{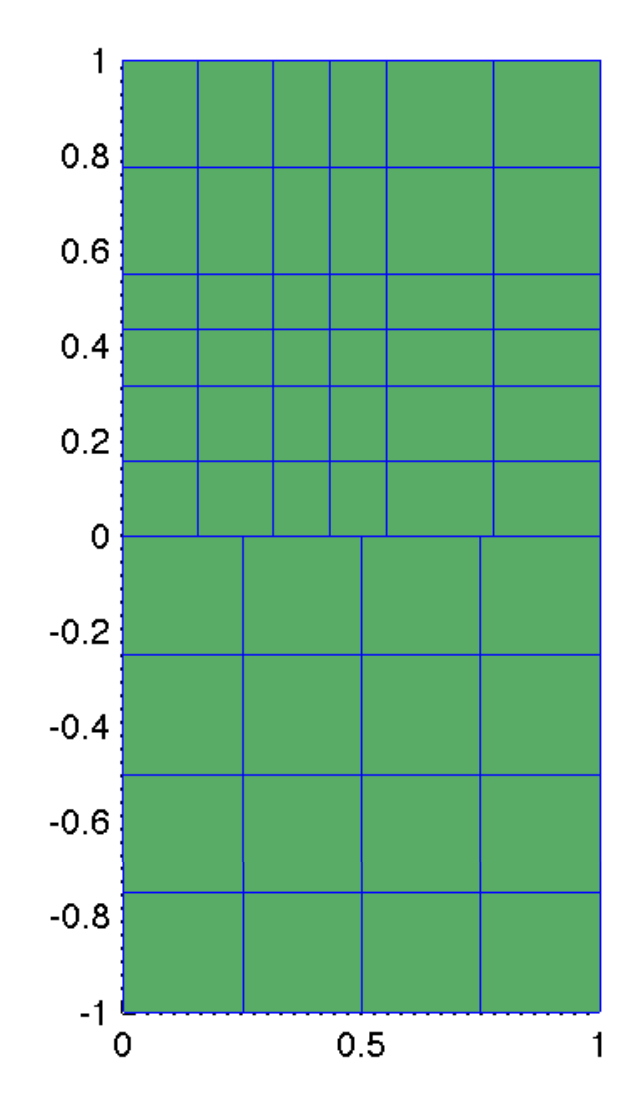}
\caption{ Different meshes at mesh refinement level $1$. From the left to the right: $M1$ to {$M3$}.}
 \label{fig:meshes}
\end{center}
\end{figure}

Regarding the meshes, we consider three different cases, presented in  Fig.~\ref{fig:meshes}.  In the first two cases, the initial master mesh is a refinement of the initial slave mesh. The initial slave mesh consists of just one element. In the case $M1$, one uniform refinement step is applied to build the master mesh, in the case $M2$ two uniform refinement steps.
Case $M3$ was chosen such that at no refinement level parts of the slave and master boundary meshes do coincide. The initial interior knots of the slave domain were chosen as $\{\pi/10, 1-\pi/7\}$ in both parametric directions, yielding $9$ elements. The initial master mesh consists of four uniform elements.

In the following, we provide different numerical error studies. We note that the inter-element smoothness of the dual functions 
 can influence the accuracy of the quadrature based on the master mesh, but not the one based on the slave mesh. Therefore for the slave integration approach, the equal order pairing with maximal smoothness is considered, i.e., $M_h = M_h^0 \subset C^{p-1}(\gamma)$, while for the non-symmetric approach we vary the dual degree. 
In all cases, the primal $L^2(\Omega)$ and the dual $L^2(\gamma)$ errors are computed by a comparison with the analytical solution stated above.

%% file: 04_1_numerics.tex
\begin{figure} 
\begin{center}
\includegraphics[width=.42\textwidth]{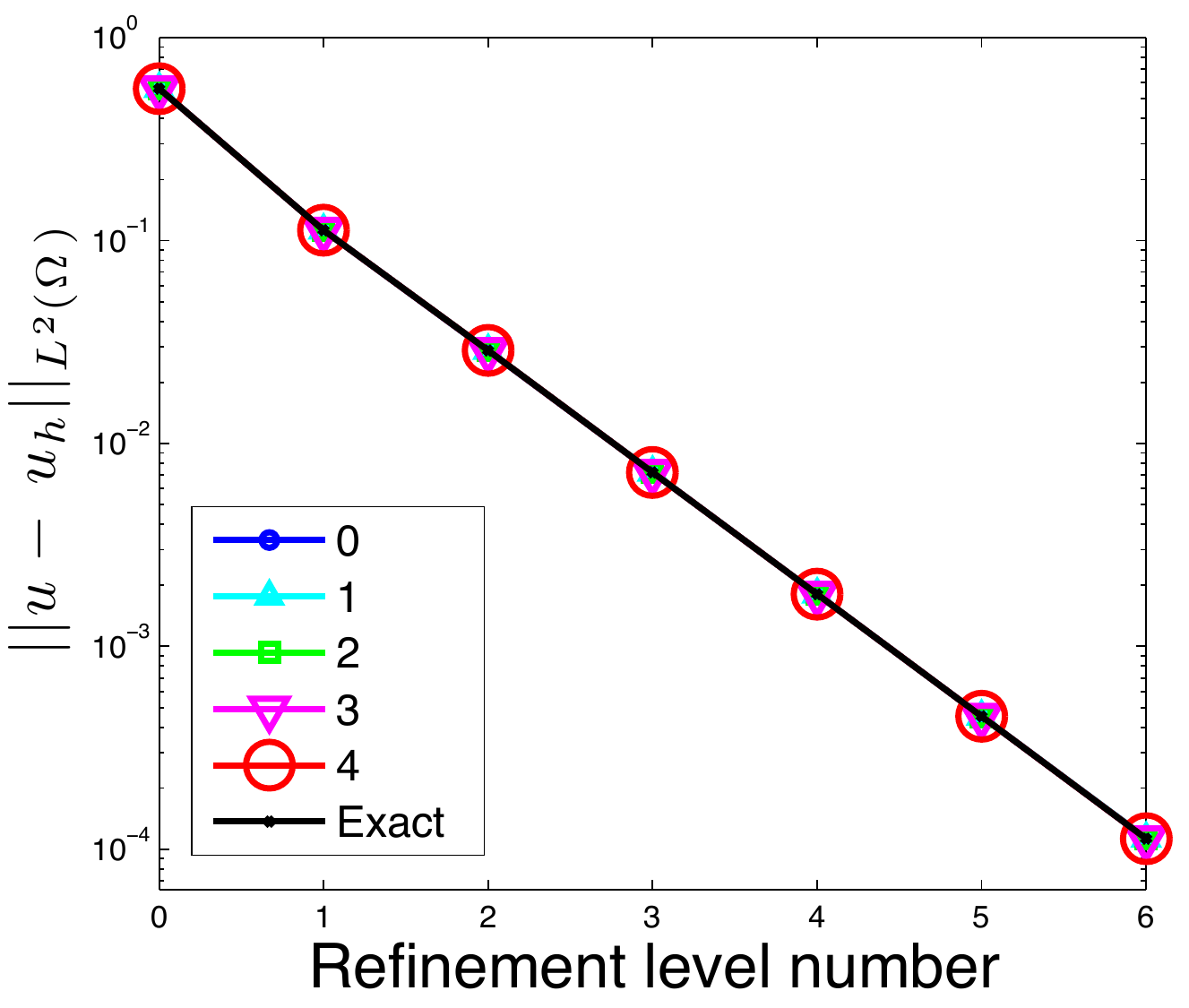}\,\,
\includegraphics[width=.42\textwidth]{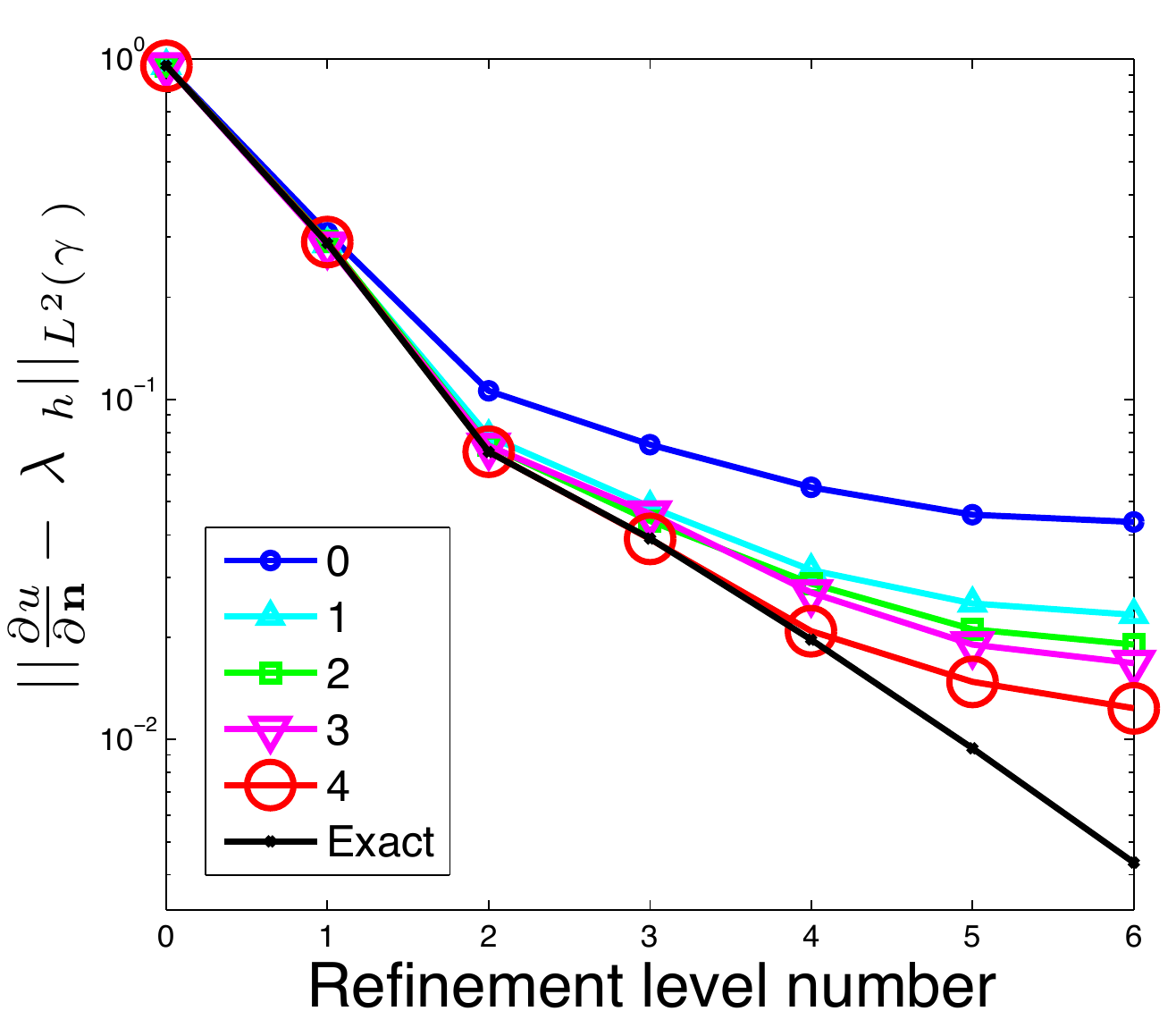}
\includegraphics[width=.42\textwidth]{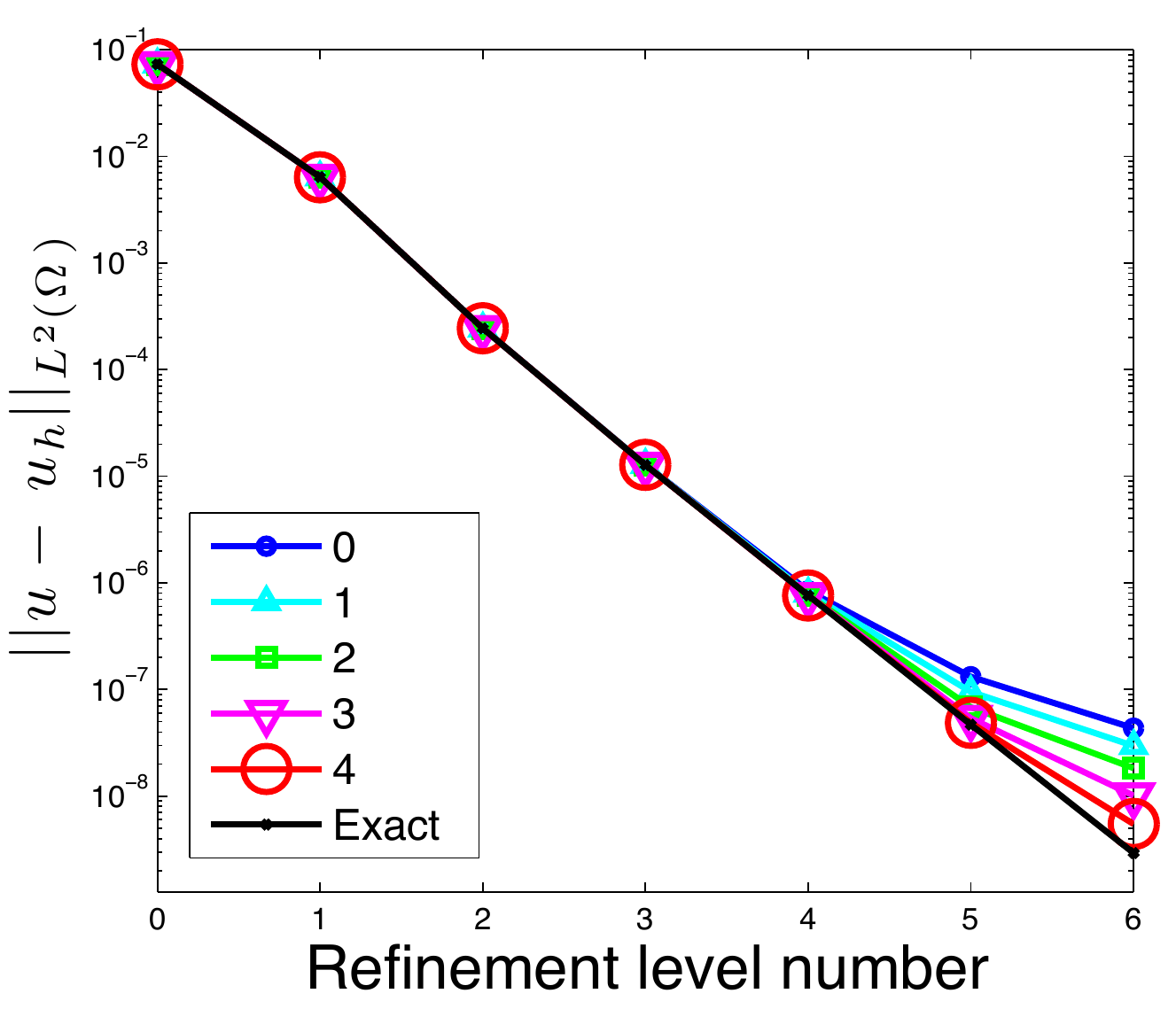}\,\,
\includegraphics[width=.42\textwidth]{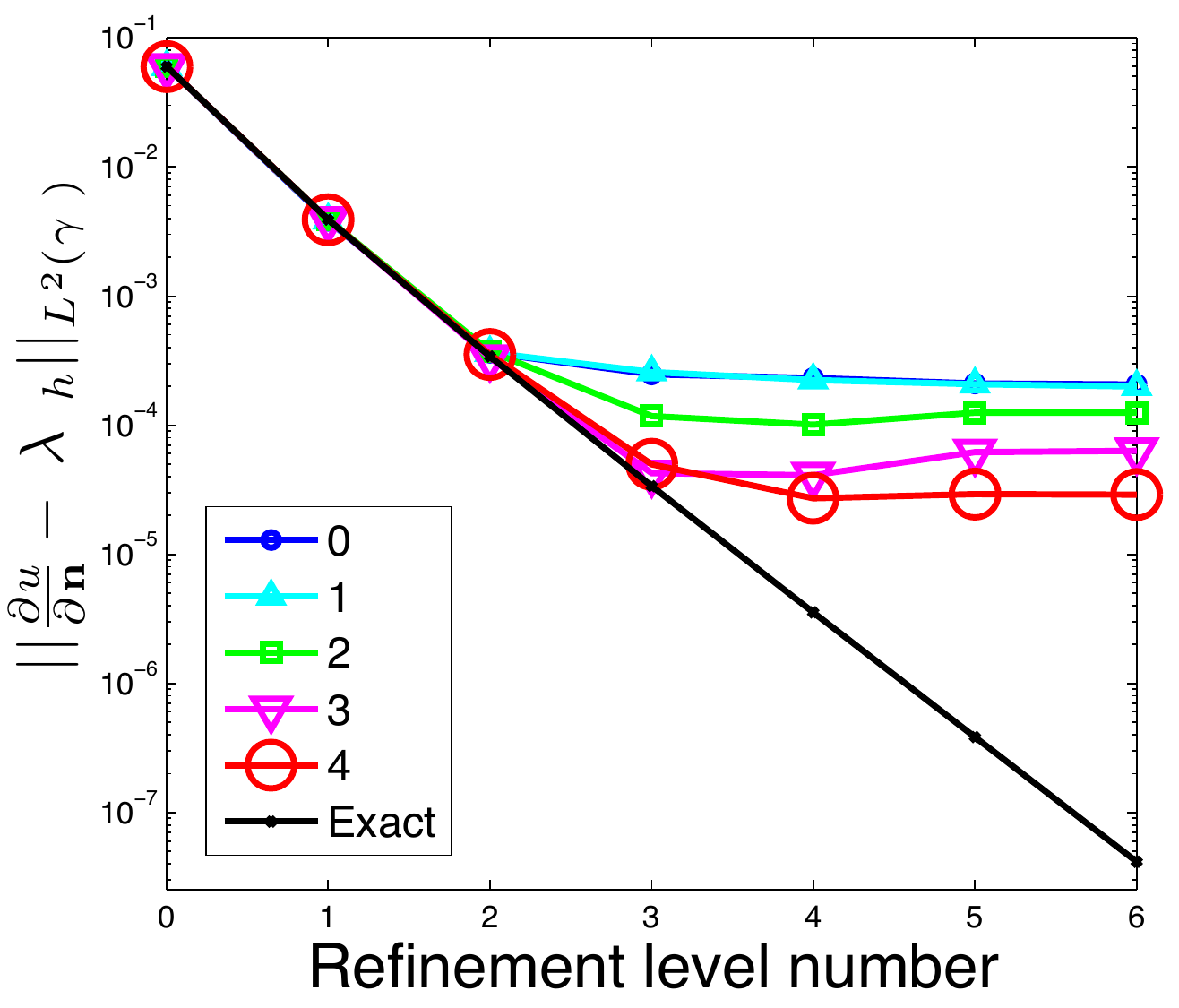}
\includegraphics[width=.42\textwidth]{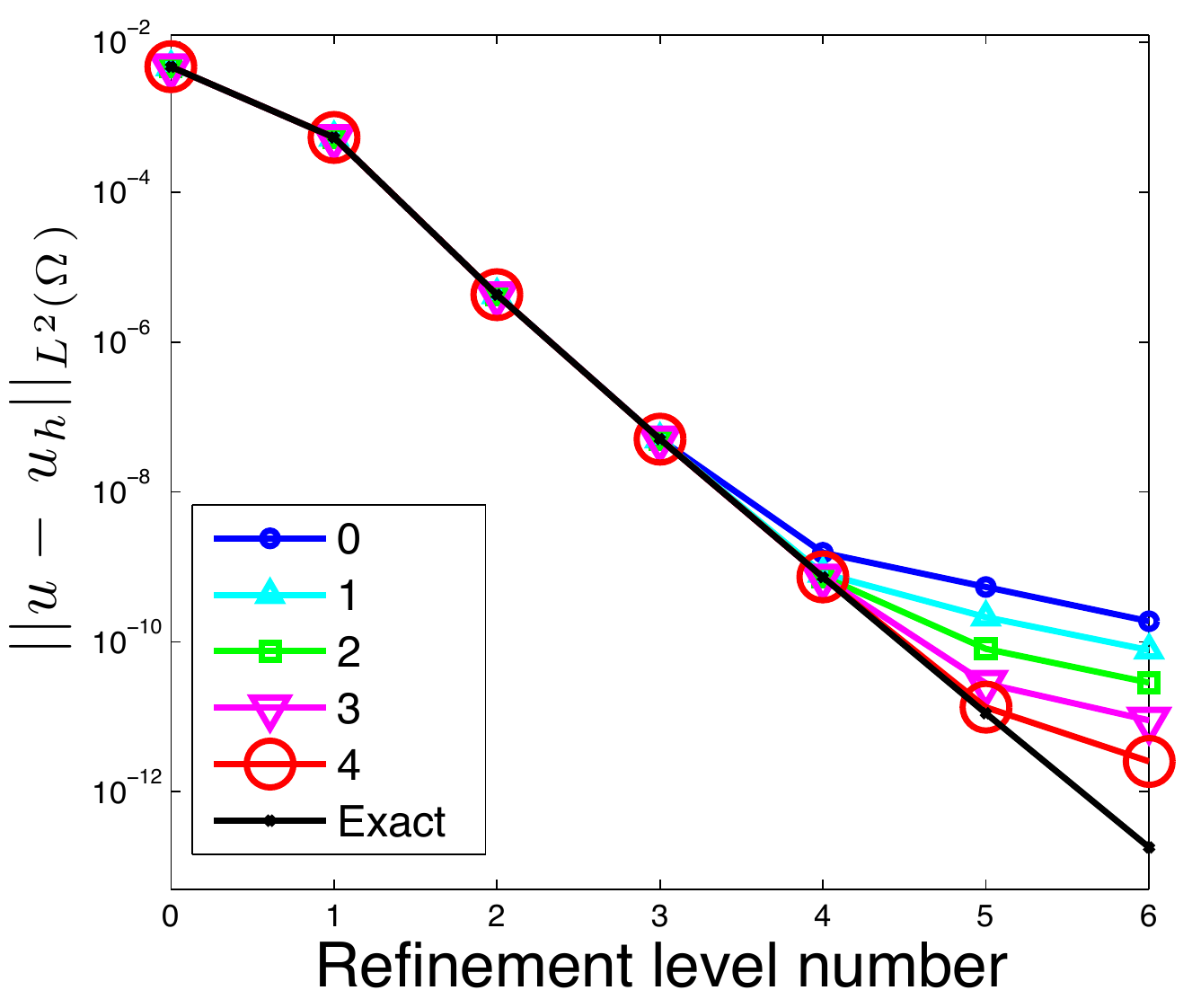}\,\,
\includegraphics[width=.42\textwidth]{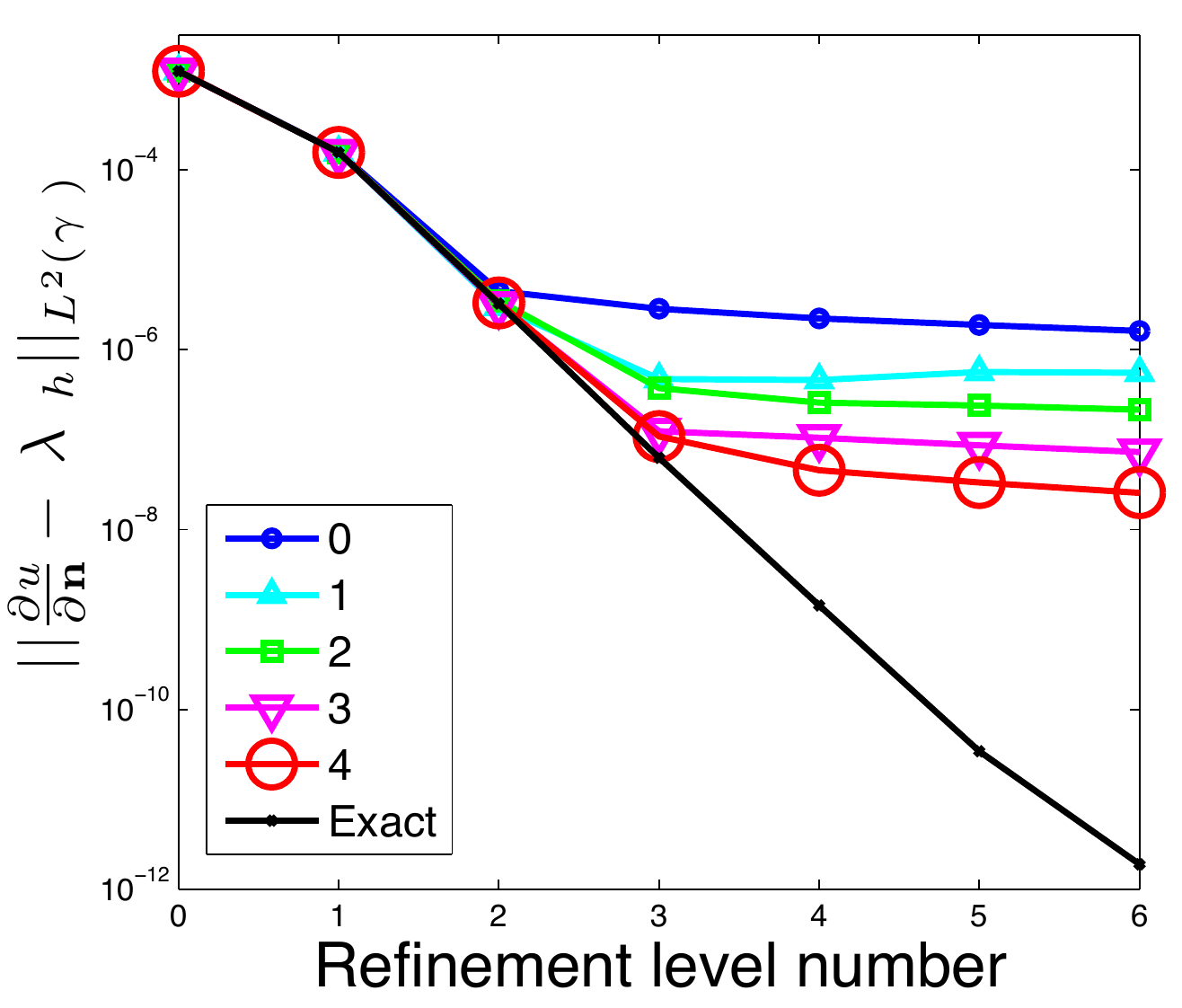}
\end{center}
\caption{2D results - $L^2$ primal (left) and dual (right) error curves for the case $M3$: { equal order} pairings { with $p=1,3,5$ (from top to bottom)} for the slave integration approach and different quadrature rule orders.}
\label{fig:diff_curves}
\end{figure}

Firstly, we consider the case $M3$, see Fig.~\ref{fig:meshes}, to measure the impact of the integration error in a general situation. A numerical error study is provided in Fig.~\ref{fig:diff_curves} for a different number of Gauss points and different spline degrees. It can clearly be seen that the primal and dual solutions are both affected by the inexact quadrature, leading to non-optimal methods. 
In all cases, the same characteristic behavior can be seen. Up to a certain refinement level, the results with inexact quadrature rules coincide with the ones with no quadrature error. Then, at a certain refinement level, the convergence order is reduced and the error is significantly larger than the exact integration one. The starting disturbance threshold is different for the primal and dual solutions as well as it differs for a different amount of quadrature points.
Moreover, in this situation the higher order splines are more disturbed by the numerical quadrature than the lower order splines. 

In almost all cases of Fig.~\ref{fig:diff_curves}, we observe poor approximation results and a reduced convergence order which is numerically independent of the spline degree.
Especially, the rate of the $L^2(\gamma)$ dual error is very low, so let us consider the final numerical convergence rate in more details. In Table~\ref{tab:eoc}, estimated convergence orders for degree $5$ and cases $M1$ and $M2$ are given. We notice that the dual $L^2(\gamma)$ rate breaks down to an order of $1/2$, while the $L^2(\Omega)$ primal rate lies about $3/2$.

\begin{table}
\begin{center}
\begin{tabular}{|l||c|c||c|c|c|}
\cline{2-5}
\multicolumn{1}{c|}{} & \multicolumn{2}{ c|| }{primal error} &  \multicolumn{2}{ c| }{dual error} \\
\cline{2-5}
\multicolumn{1}{c|}{} & \multicolumn{1}{c|}{case $M1$}  & \multicolumn{1}{c||}{case $M2$}   & \multicolumn{1}{c|}{case $M1$}  & \multicolumn{1}{c|}{case $M2$} 
 \\ \hline
quad. rule order 0& 1.63 & 1.74  & 0.50 & 0.50   \\
quad. rule order 1& 1.63 & 1.54  & 0.50 & 0.50 \\
quad. rule order 2& 1.63 & 1.55  & 0.50 & 0.50 \\
quad. rule order 3& 1.63 & 1.58  & 0.50 & 0.50 \\
quad. rule order 4& 1.63 & 1.56  & 0.50 & 0.50 \\
quad. rule order 5& 1.63 & 1.50  & 0.50 & 0.50 \\
\hline
\end{tabular}
\caption{2D results - Last estimated order of convergence of the primal and dual $L^2$ errors for the cases $M1$ and $M2$: pairing $P5-P5$ for the slave integration approach and different quadrature rule orders.  }
\label{tab:eoc}
\end{center}
\end{table}

\begin{figure} 
\begin{center}
\includegraphics[width=.45\textwidth]{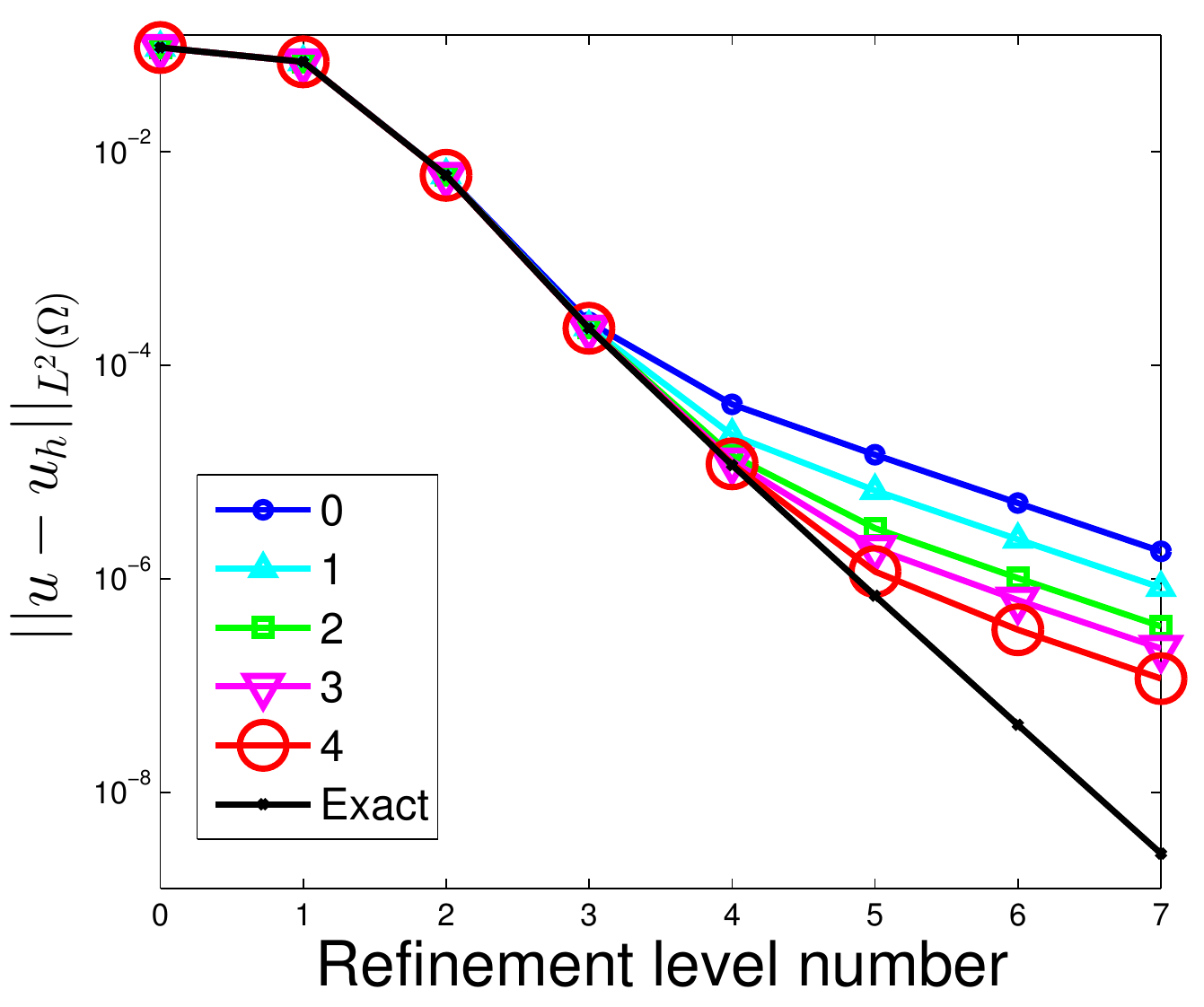}\,\,
\includegraphics[width=.45\textwidth]{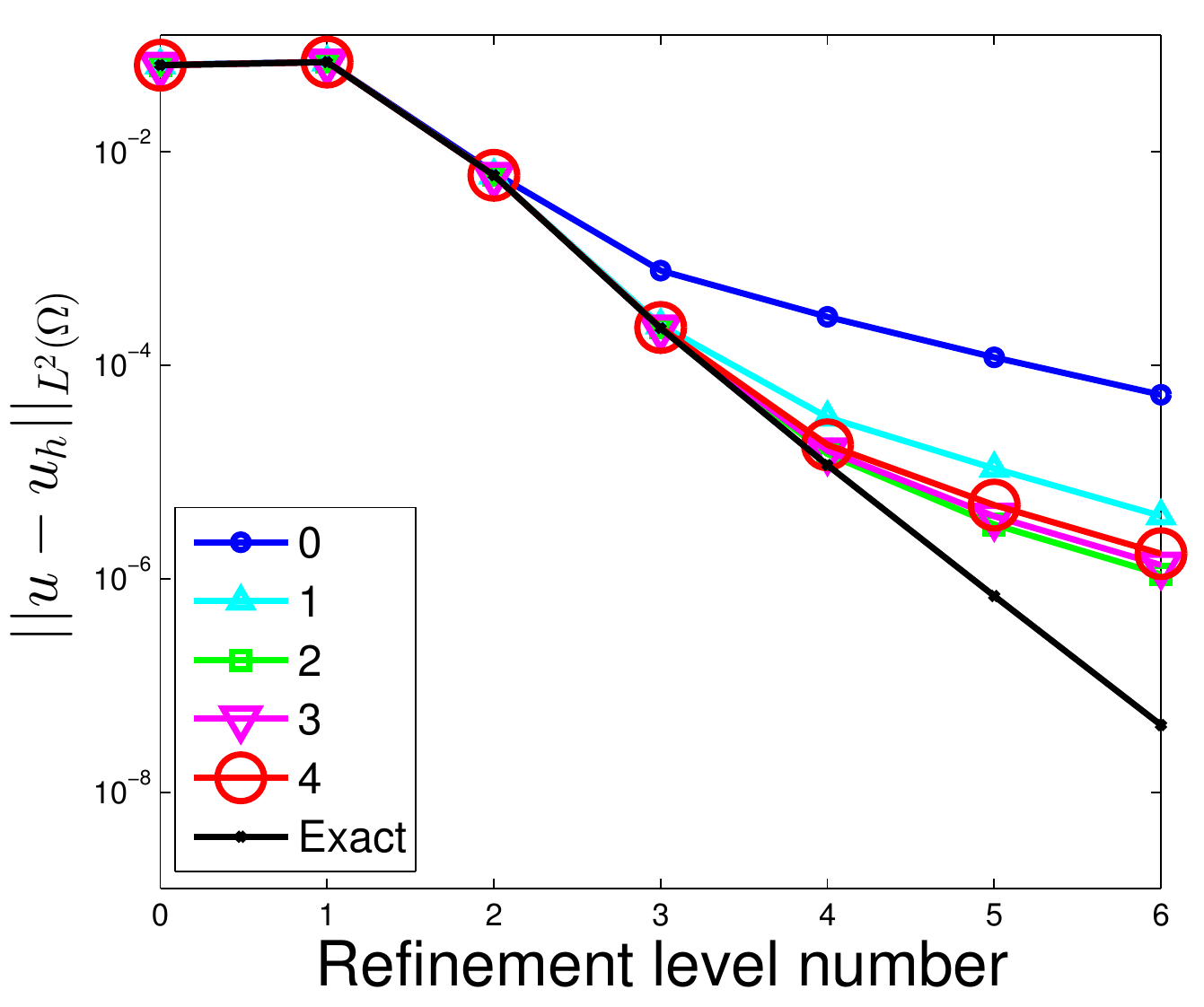}
\includegraphics[width=.45\textwidth]{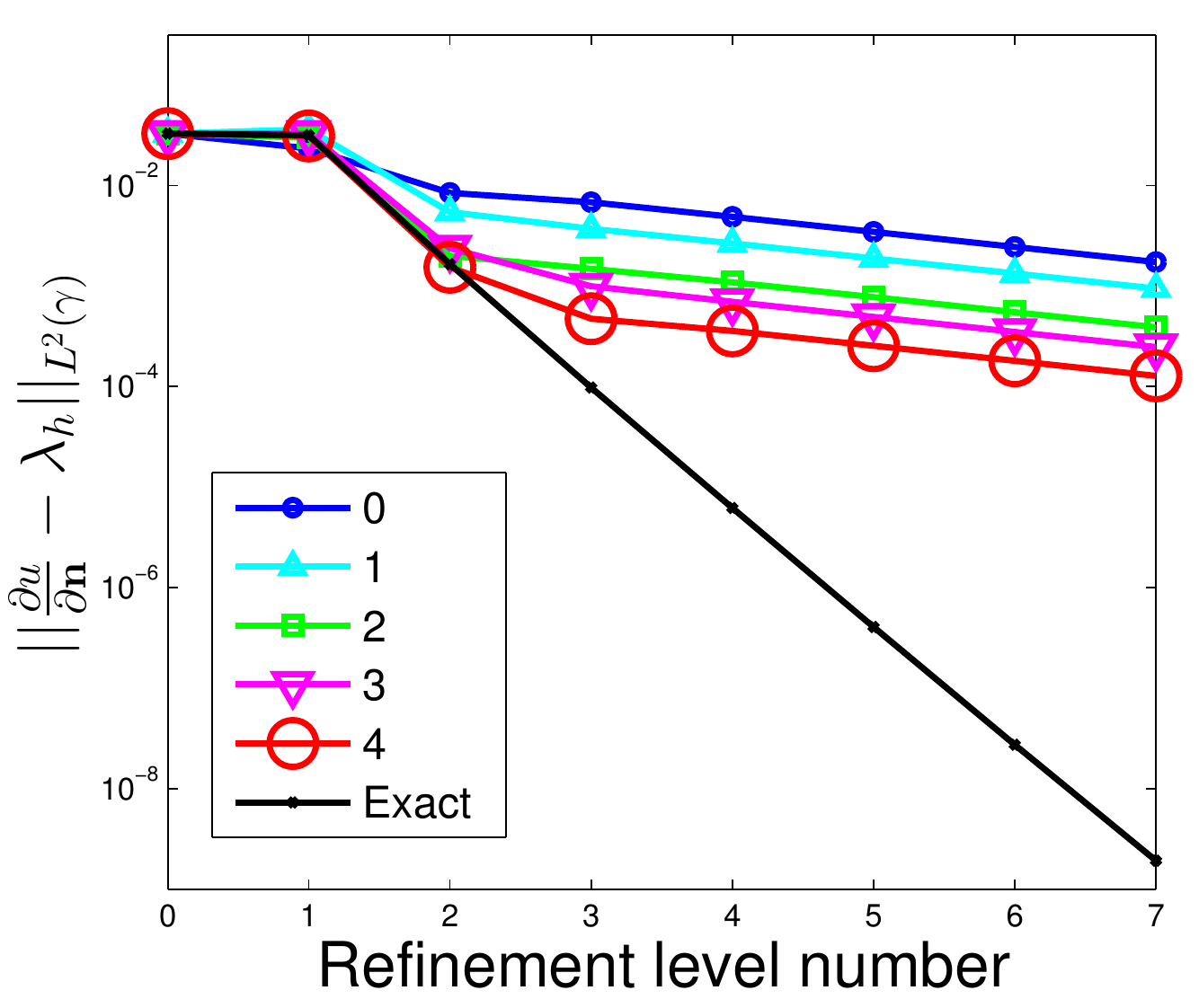}\,\,
\includegraphics[width=.45\textwidth]{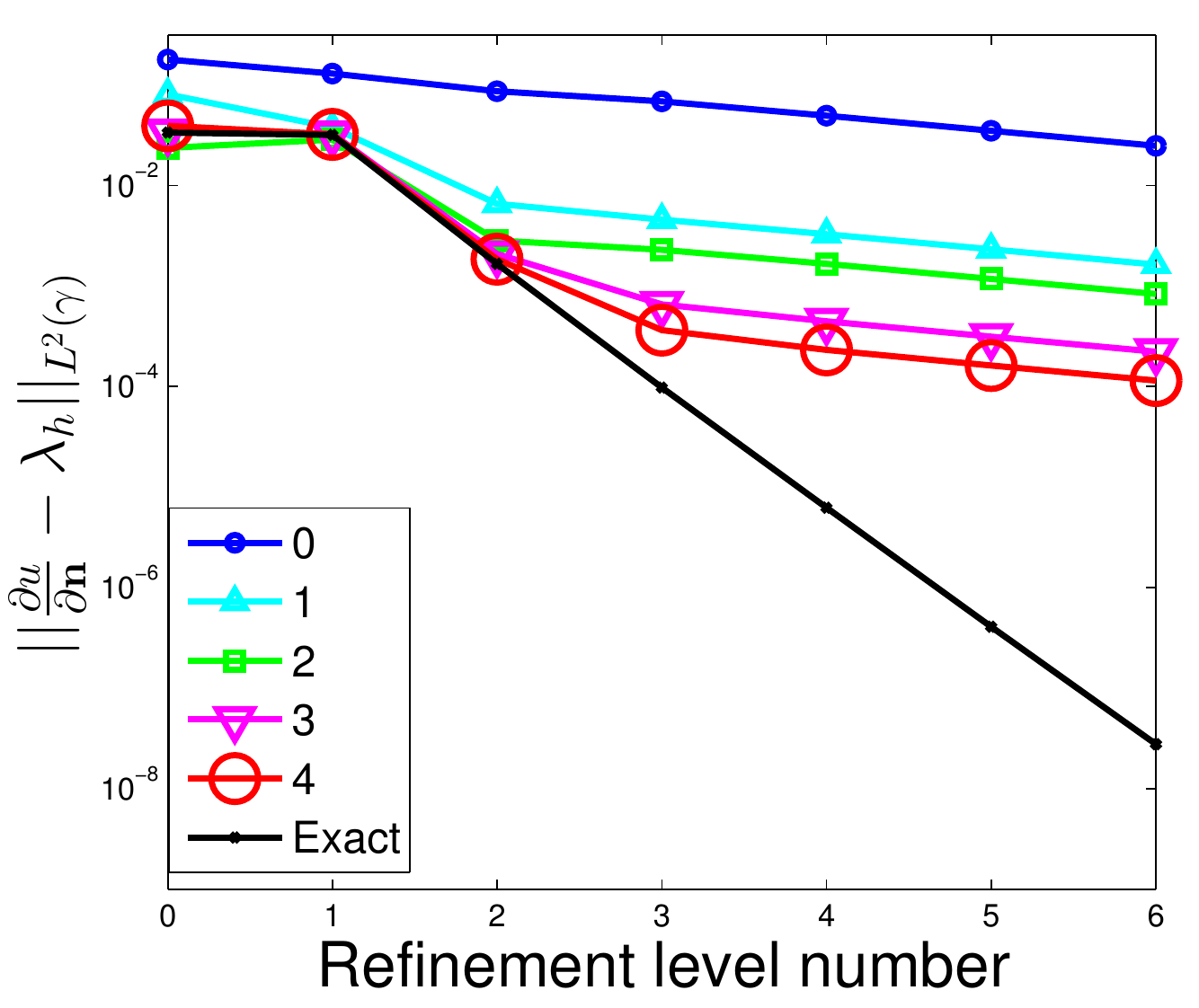}
\end{center}
\caption{2D results - $L^2$ primal (top) and dual (bottom) error curves for the cases $M1$ (left) and $M2$ (right): pairing $P3-P3$ for the slave integration approach and different quadrature rule orders.}
\label{fig:half_quarter_curves}
\end{figure}

Secondly, we consider an even more simple situation to  show that even then the impact of the slave integration is noticeable. Let us focus on the cases $M1$ and $M2$, see Fig.~\ref{fig:meshes}, for which the master mesh is a refinement of the slave mesh. 
See Fig.~\ref{fig:half_quarter_curves} for a comparison of results between the cases $M1$ and $M2$ for a spline degree $p=3$. 
We note that the low convergence order of the primal and dual solutions, as remarked above, already appear in this simple context. 
{
Moreover, for a fixed number of slave elements, the error is increasing with the number of master elements.  This is expected as there are more points of reduced smoothness which are not taken account by the quadrature rule.
}

{
Thirdly, we have additionally compared the case $M3$ with a similar situation in which the master and slave roles are inverted. The results also show that the integration error is increasing with the
increase of the master element number. Thus, in accordance to the practical applications, in a slave integration context it seems worthwhile to choose the slave domain as the finest one.
}

\begin{figure} 
\begin{center}
\includegraphics[width=.42\textwidth]{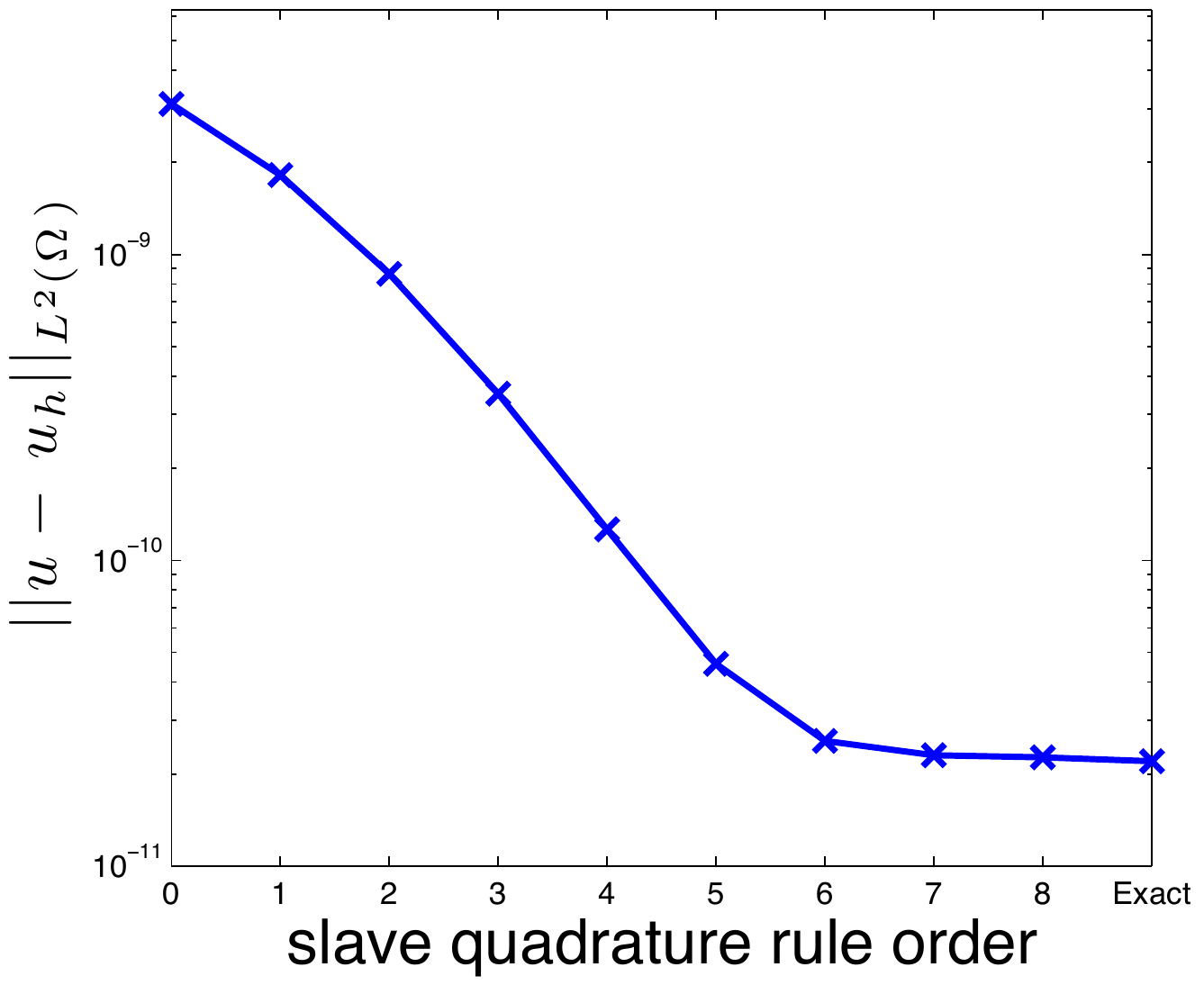}\,\,
\includegraphics[width=.42\textwidth]{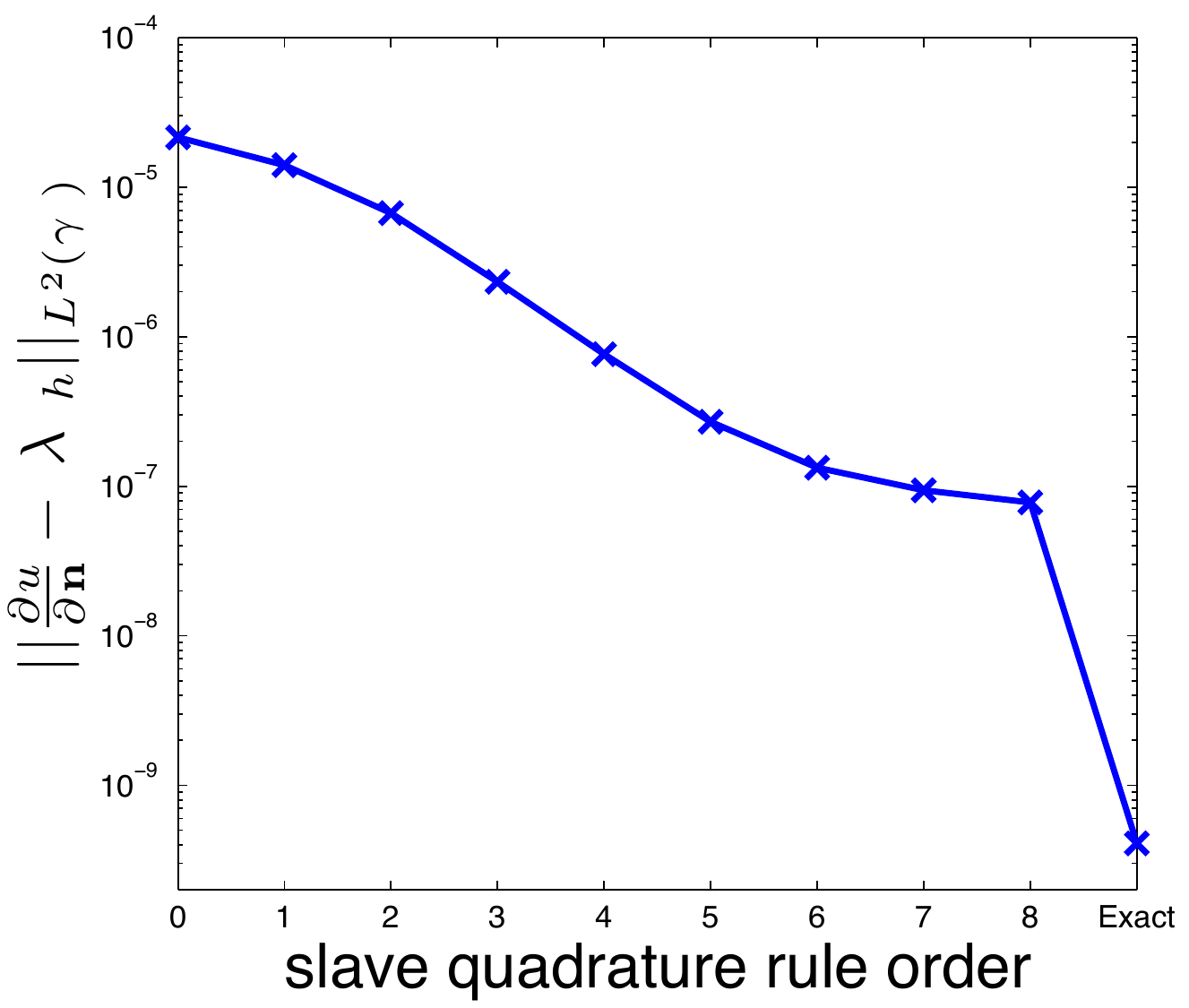}
\end{center}
\caption{2D results - $L^2$ primal (left) and dual (right) errors at refinement level number 6 as a function of the quadrature rule order for the case $M3$: pairing $P3-P3$ for the slave integration approach.}
\label{fig:diff_quad_curves}
\end{figure}

Moreover, it can be observed that on coarse meshes using the slave integration method it is possible to recover  the accuracy of the optimal mortar method simply by increasing the number of quadrature points, see Fig.~\ref{fig:diff_quad_curves}. 
 However, it has also been shown that the number of necessary quadrature points is drastically increasing with the refinement level. It can easily be seen that the number of Gauss points gets soon impracticably large, see the right picture of Fig.~\ref{fig:diff_quad_curves}. Furthermore, in several cases, the disturbance to the mortar method has been observed to be more severe for higher order functions.

%% file: 04_2_numerics.tex
The non-symmetric saddle point problem based on the two different quadrature rules, see Sec.~\ref{sec:3}, was introduced to overcome the non-optimality of the pure slave integration approach in a finite element context. Due to the suboptimal results seen in the previous section, it is also interesting to consider it in an isogeometric context.

\begin{figure}
\begin{center}
\includegraphics[width=.42\textwidth]{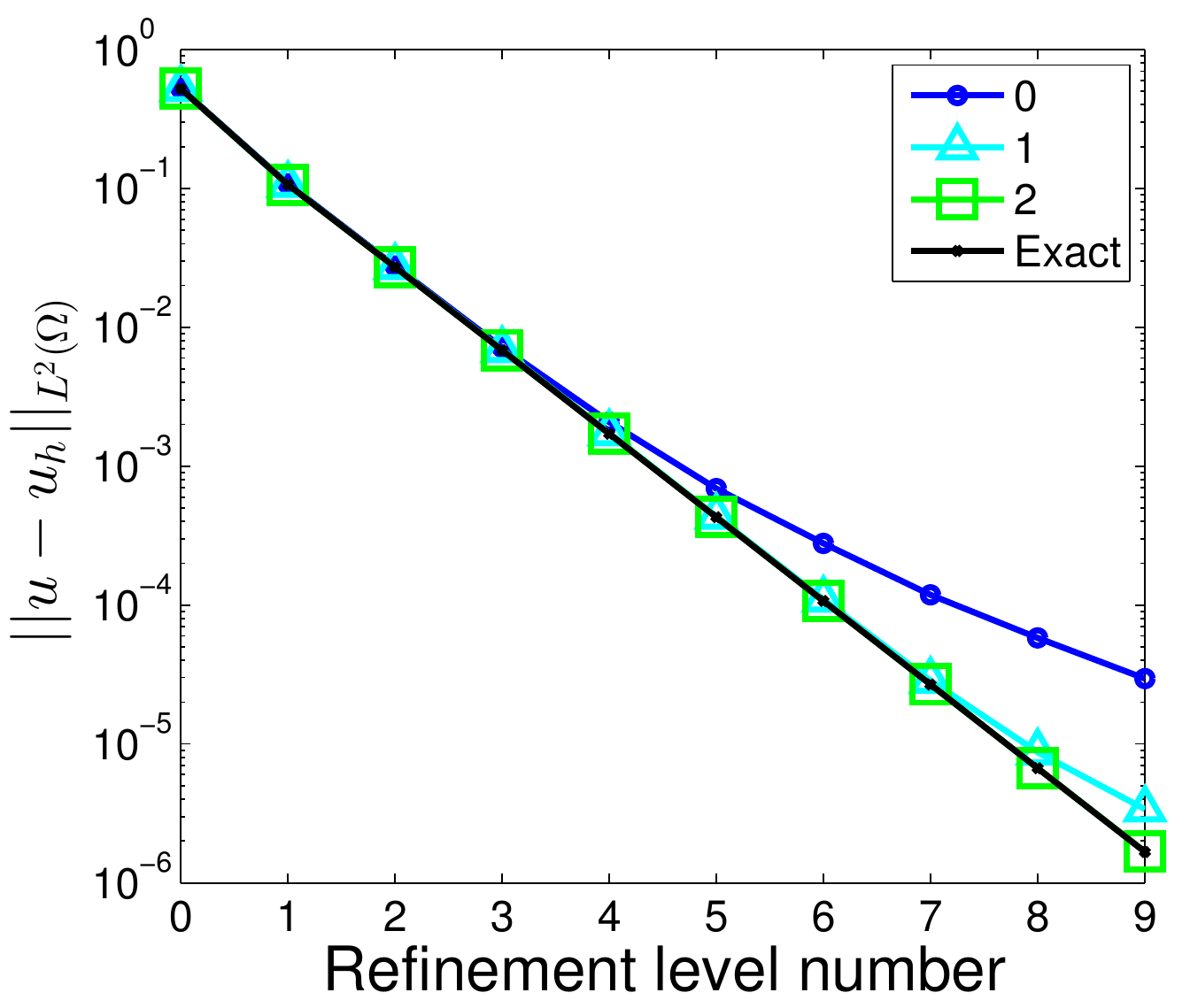}\,\,
\includegraphics[width=.42\textwidth]{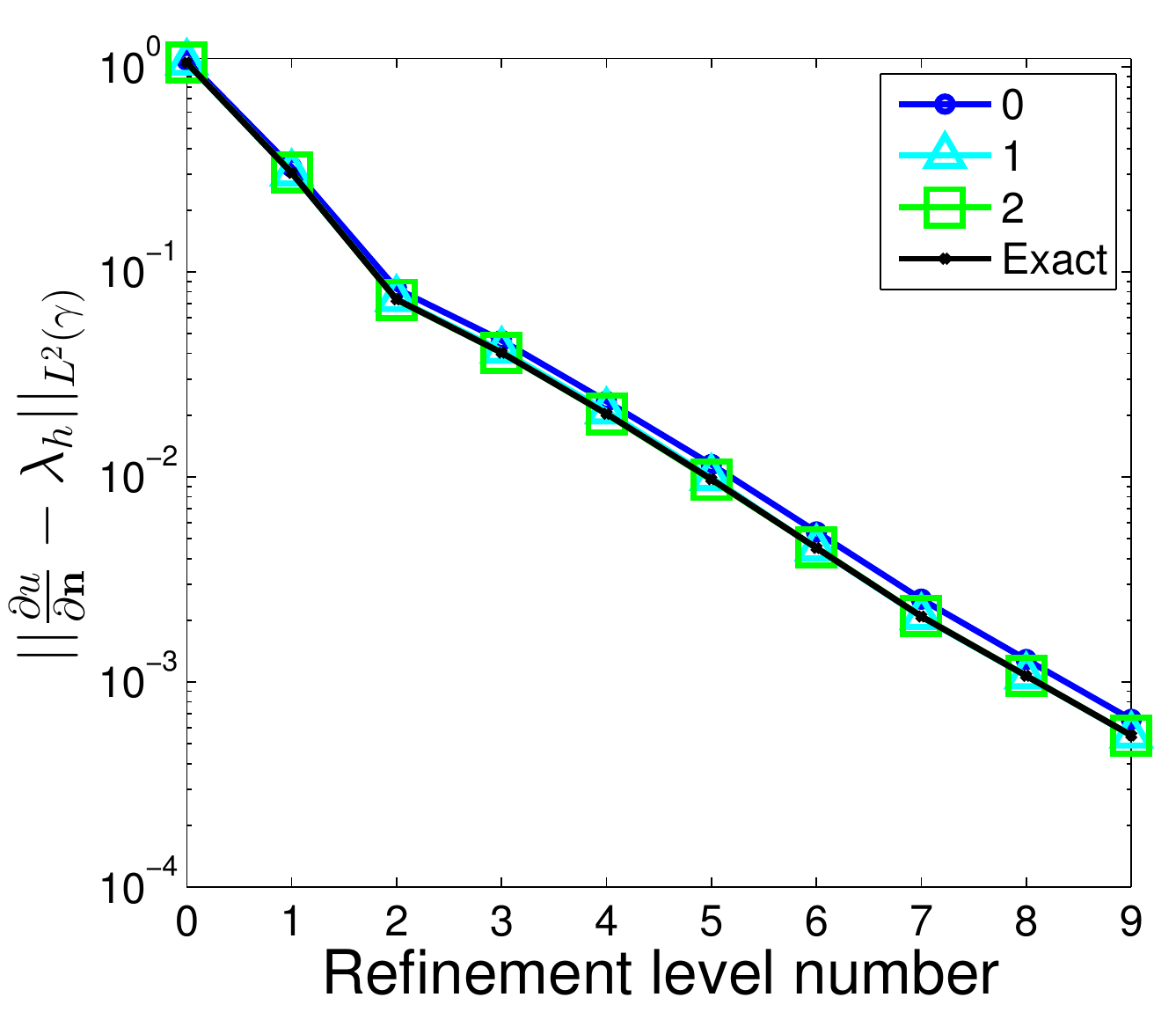}
\end{center}
\caption{2D results - $L^2$ primal (left) and dual (right) error curves for the case $M3$: equal order pairing $p=1$ for the non-symmetric approach and different quadrature rule orders.}\label{fig:pg_p1}
\end{figure}

Firstly, we consider same degree pairings. In almost all tested cases the results of the non-symmetric approach are comparable to the results of the exact integration case. However, we note that differences could still be seen in some cases. For example, for a degree $p=1$ in the case $M3$, we obtained a non-optimal method, see in Fig.~\ref{fig:pg_p1} the corresponding primal and dual error curves. Note that we do not show any curves in the cases where no disturbance is observed. For example for degree $p=5$ we observed convergence almost up to machine precision without any remarkable difference compared to the exact integration case.

\begin{figure}
\begin{center}
\includegraphics[width=.42\textwidth]{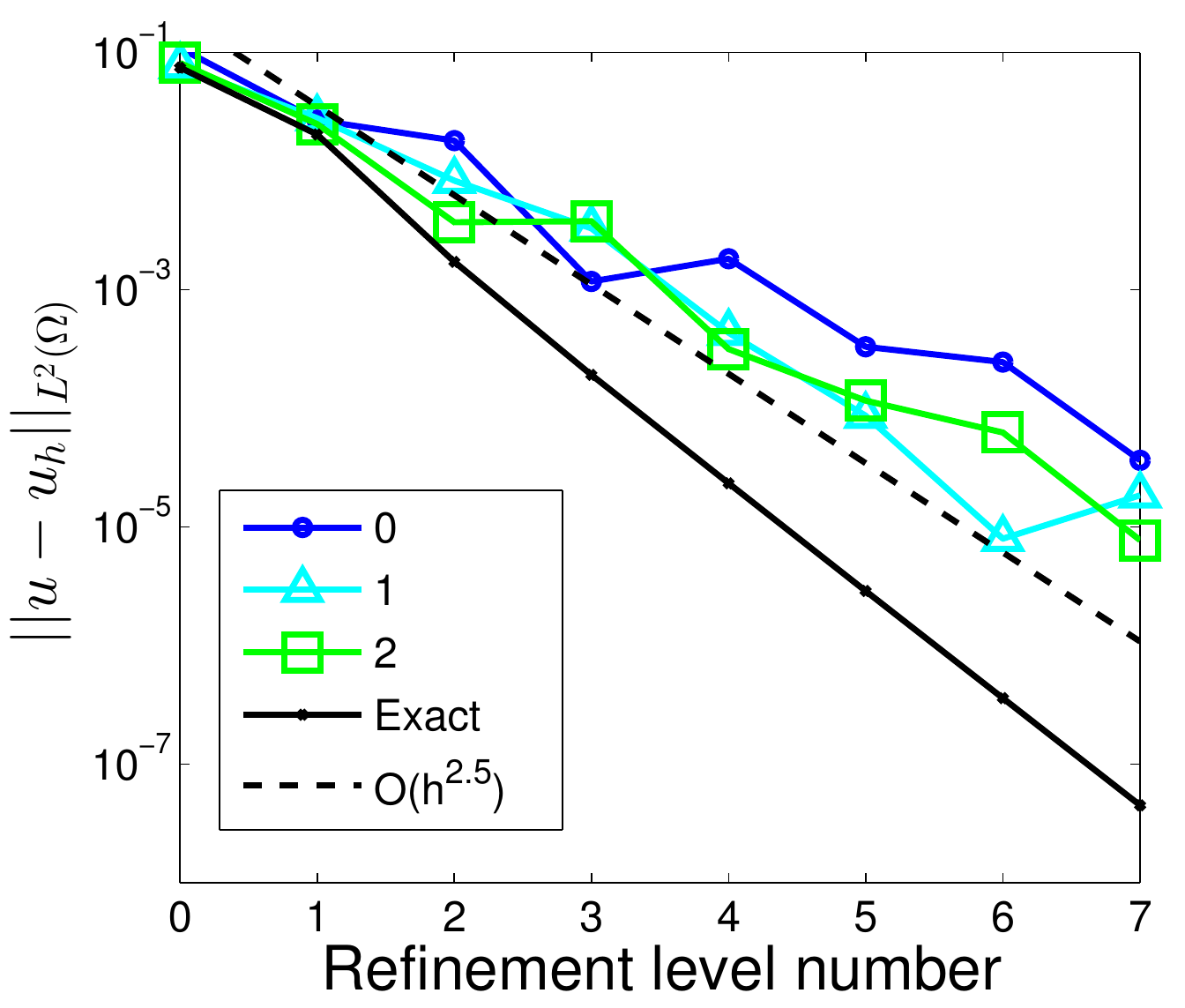}\,\,
\includegraphics[width=.42\textwidth]{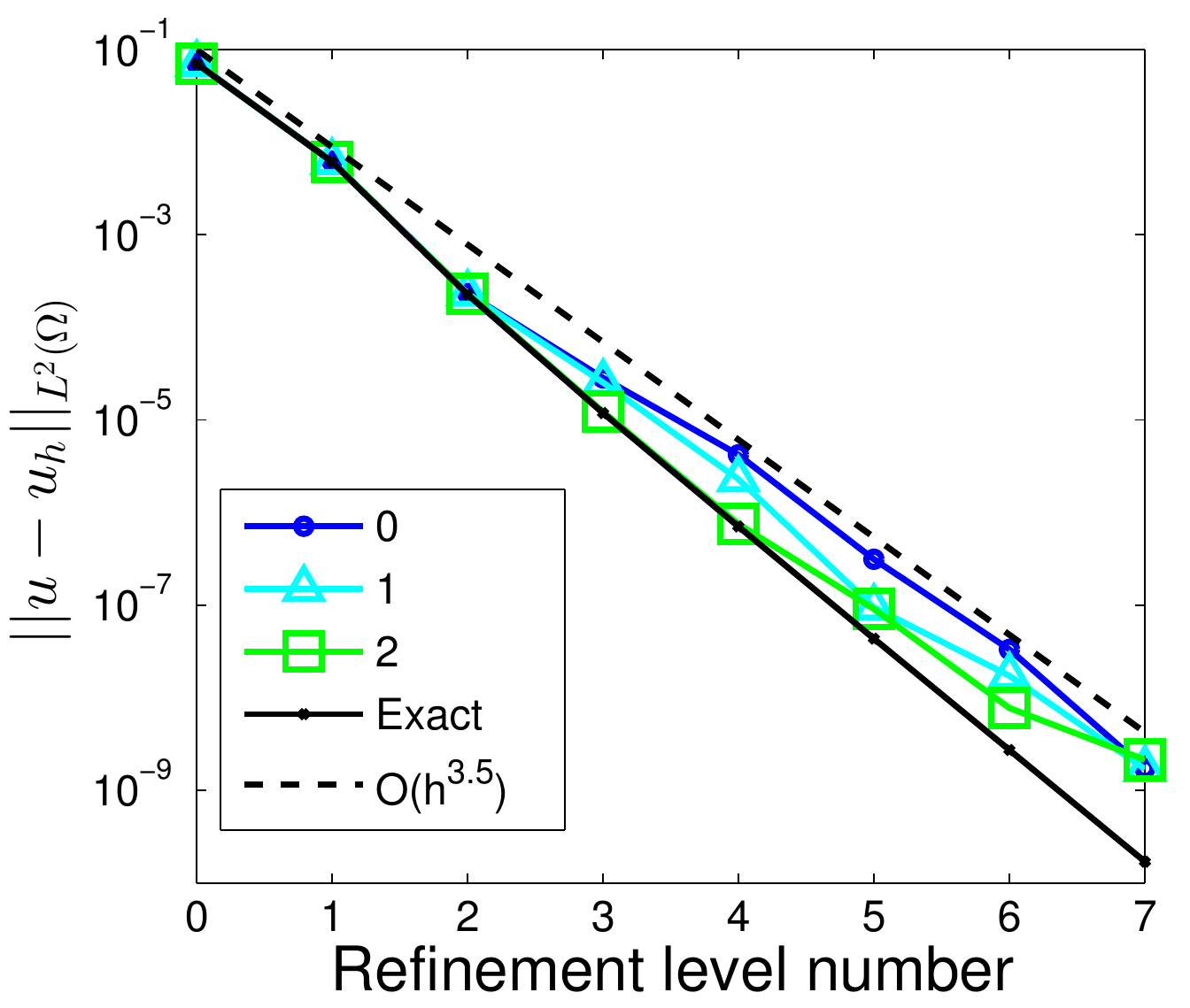}\\
\includegraphics[width=.42\textwidth]{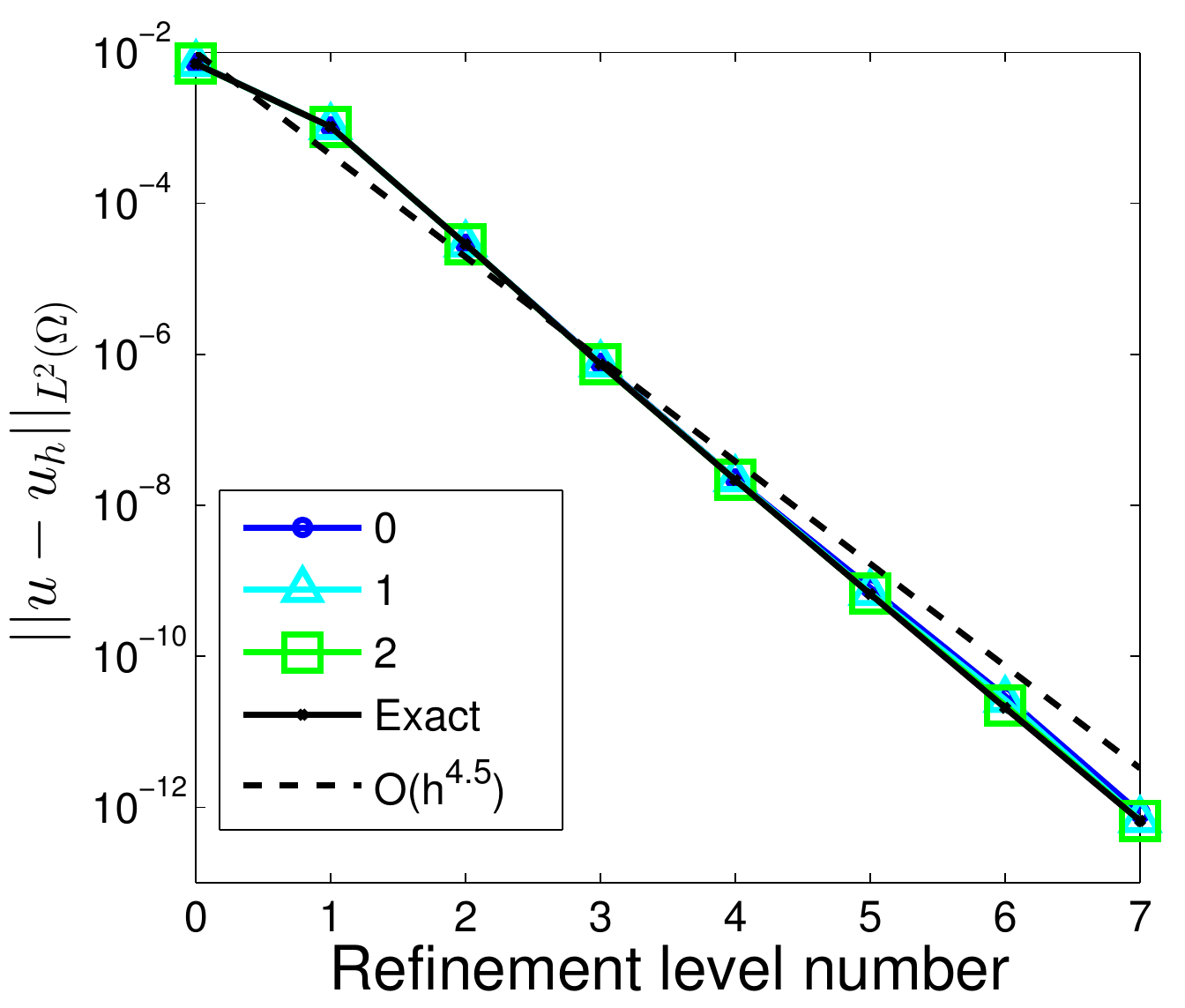}\,\,
\includegraphics[width=.42\textwidth]{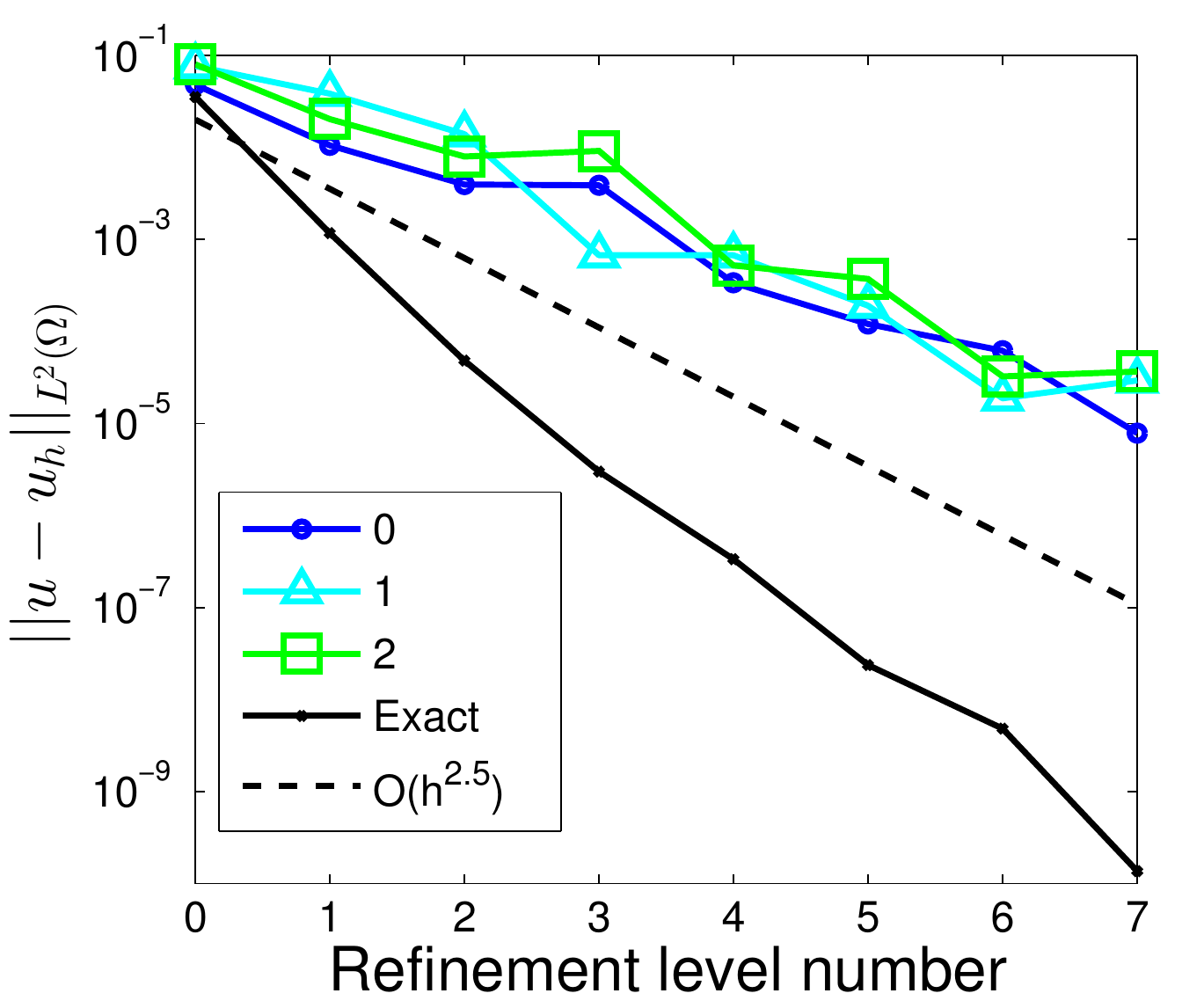}
\end{center}
\caption{2D results - $L^2$ primal error curves for the case $M3$: different order pairings for the non-symmetric approach and different quadrature rule orders. Top left: $P2-P0$. Top right: $P3-P1$. Bottom left: $P4-P2$. Bottom right: $P4-P0$.}\label{fig:pg_different_degrees}
\end{figure}

Secondly, we consider dual spaces with lower degrees than the primal ones. Note that in~\cite{brivadis:14} stability for these pairings was only observed if the primal and the dual degrees have the same parity.
Similar to the equal order case, the dual error is not affected. In Fig.~\ref{fig:pg_different_degrees} primal error curves are shown for all stable different degree pairings up to a primal degree $p=4$. We note that theoretically, we expect sub-optimal primal error rates even in the exact integration case, although often improved convergence rates were observed. For a dual degree $p-2k, \, k \in \mathbb{N}$, we can expect a convergence of order $\mathcal{O} (h^{p-2k+5/2})$ in the $L^2(\Omega)$ norm, see the dashed lines in Fig.~\ref{fig:pg_different_degrees}.
For the $P4-P2$ and $P3-P1$ pairings, we observe slight differences compared to the exact integration results, but note that the convergence rate is not far from the theoretical expectation. The situation is different for the $P4-P0$ and $P2-P0$ pairings, for which the rate is more disturbed and even below the theoretical expectation. This could be explained by the discontinuity of the dual basis functions, which introduces large errors in the integration approximation done with a rule based on the master mesh, which does not respect these discontinuities.

%% file: 04_3_numerics.tex
\begin{figure}\centering
\begin{center}
\input{figures/axe_picture}
\includegraphics[width=.4\textwidth]{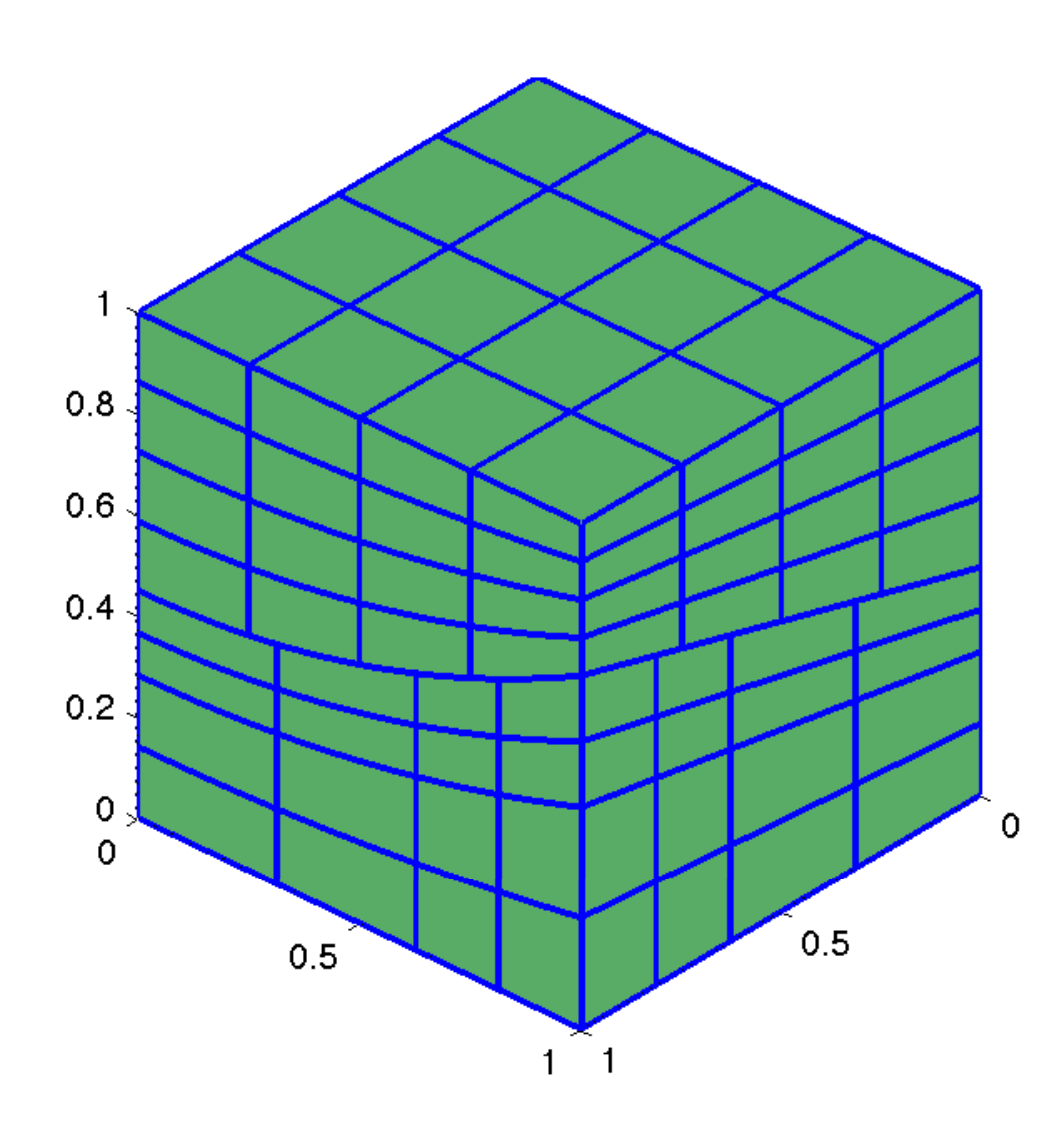}\,
\includegraphics[width=.3\textwidth]{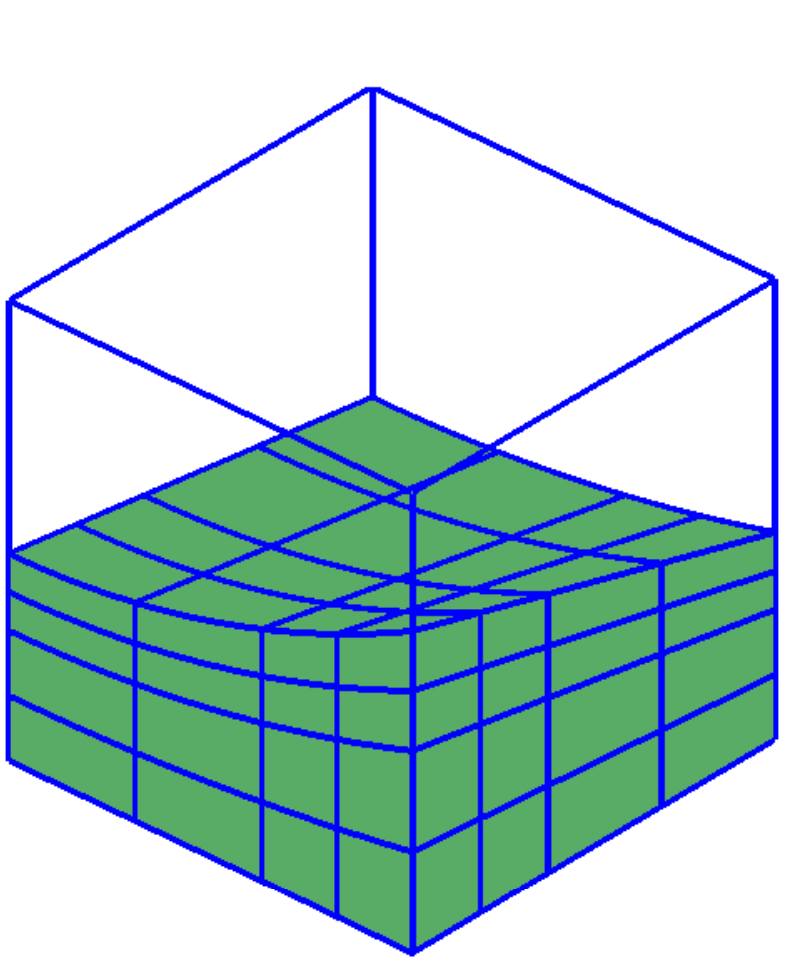}
\end{center}
\caption{Meshes at refinement level 1 (left) and the slave domain (right) illustrating the curved interface.} \label{fig:3d_geometry}
\end{figure}

As a second example, we consider a three-dimensional problem with a curved interface. Precisely, we consider the Poisson problem $-\Delta u = f$ on the domain $\Omega = (0,1)^3$, which is divided into two patches by the interface $\gamma = \{ (x,y,\rho(x)), (x,y) \in (0,1)^2
\}$, with $\rho(x,y) =1/8 \,(1+x)(1+y^2)+1/5$, see Fig.~\ref{fig:3d_geometry}. The bottom domain is set as the slave domain. The internal load and the boundary conditions are manufactured to have for analytical solution
\[
u(x,y,z) = \cos(2\pi x)\cos(2\pi y)\sin(2\pi z).
\]
Note that due to the curved interface, the normal derivative has a complicated form, but is still explicitly computable. Neumann conditions are applied such that no cross point modification is necessary.
The initial master mesh has 8 uniform elements, while the initial slave mesh has 8 elements given by the breakpoint vector $\{0, \pi/5, 1\}$ in each direction. In the following, we provide some numerical error studies, considering the slave integration approach as well as the non-symmetric approach.

\begin{figure}[htbp]
\begin{center}
\includegraphics[width=.42\textwidth]{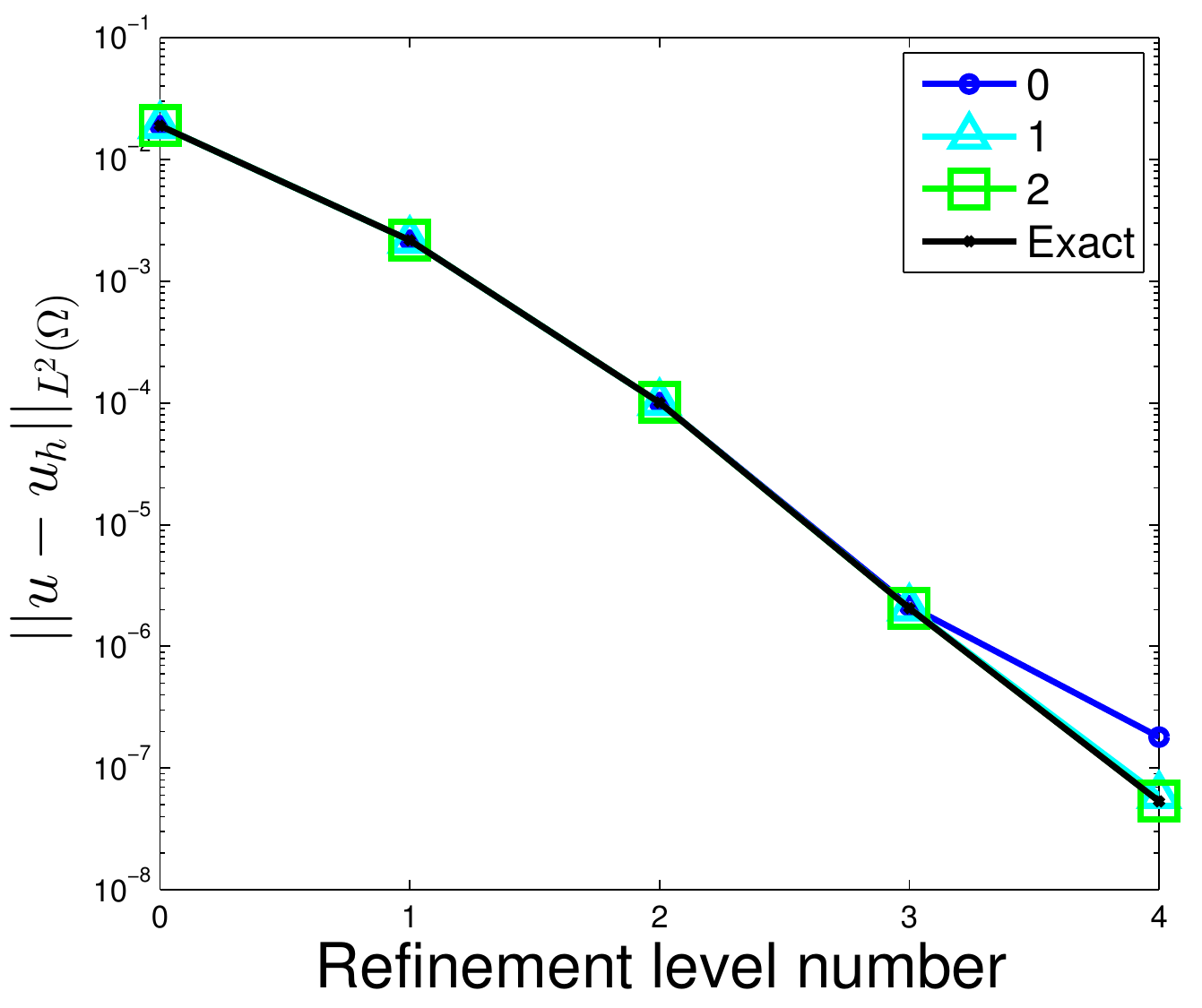}\,\,
\includegraphics[width=.42\textwidth]{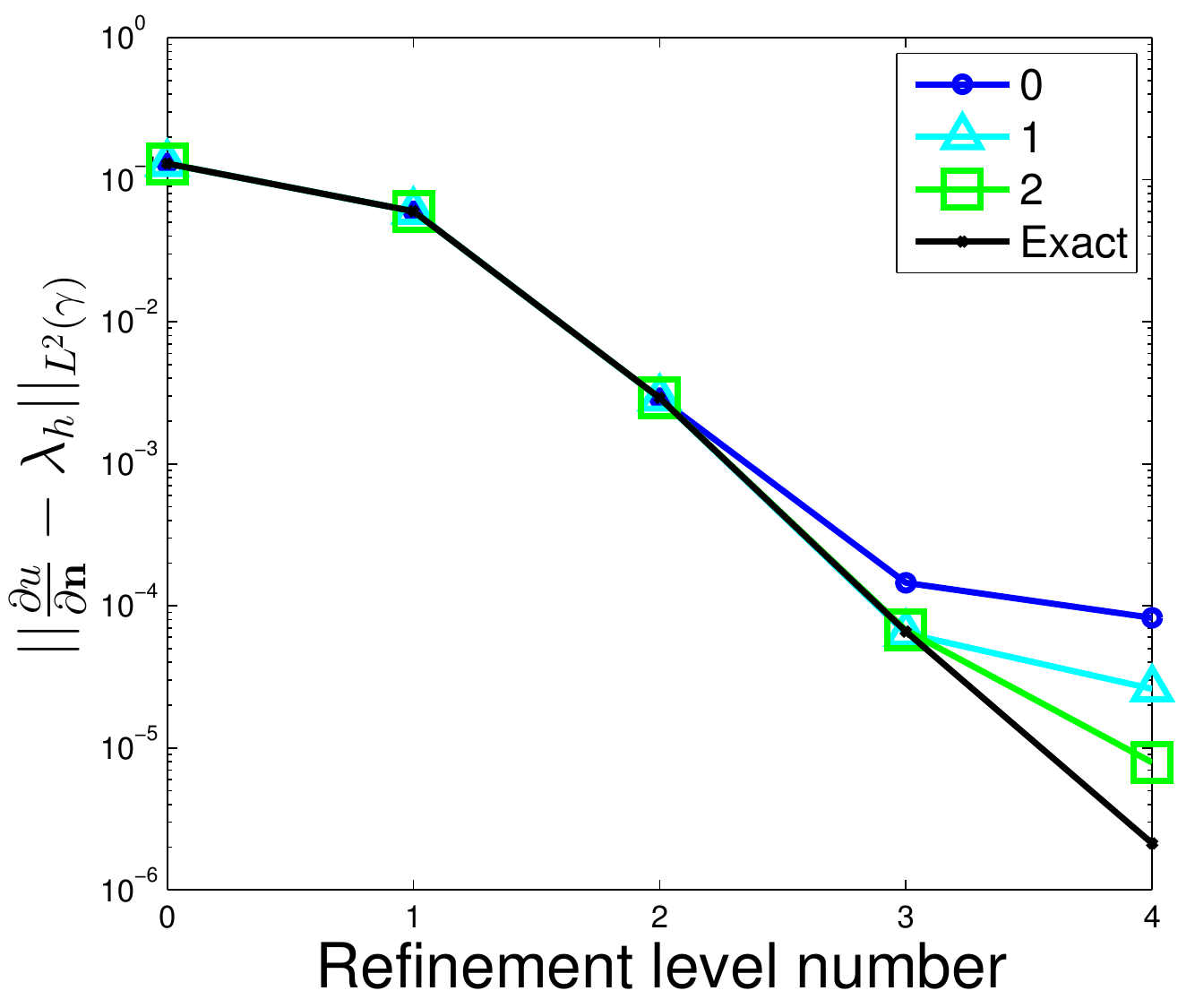}
\end{center}
\caption{3D results - $L^2$ primal (left) and dual (right) error curves for the pairing $P4-P4$, for the slave integration approach and different quadrature rule orders.}\label{fig:3d_slave_p4}
\end{figure}

\begin{figure}[htbp]
\begin{center}
\includegraphics[width=.42\textwidth]{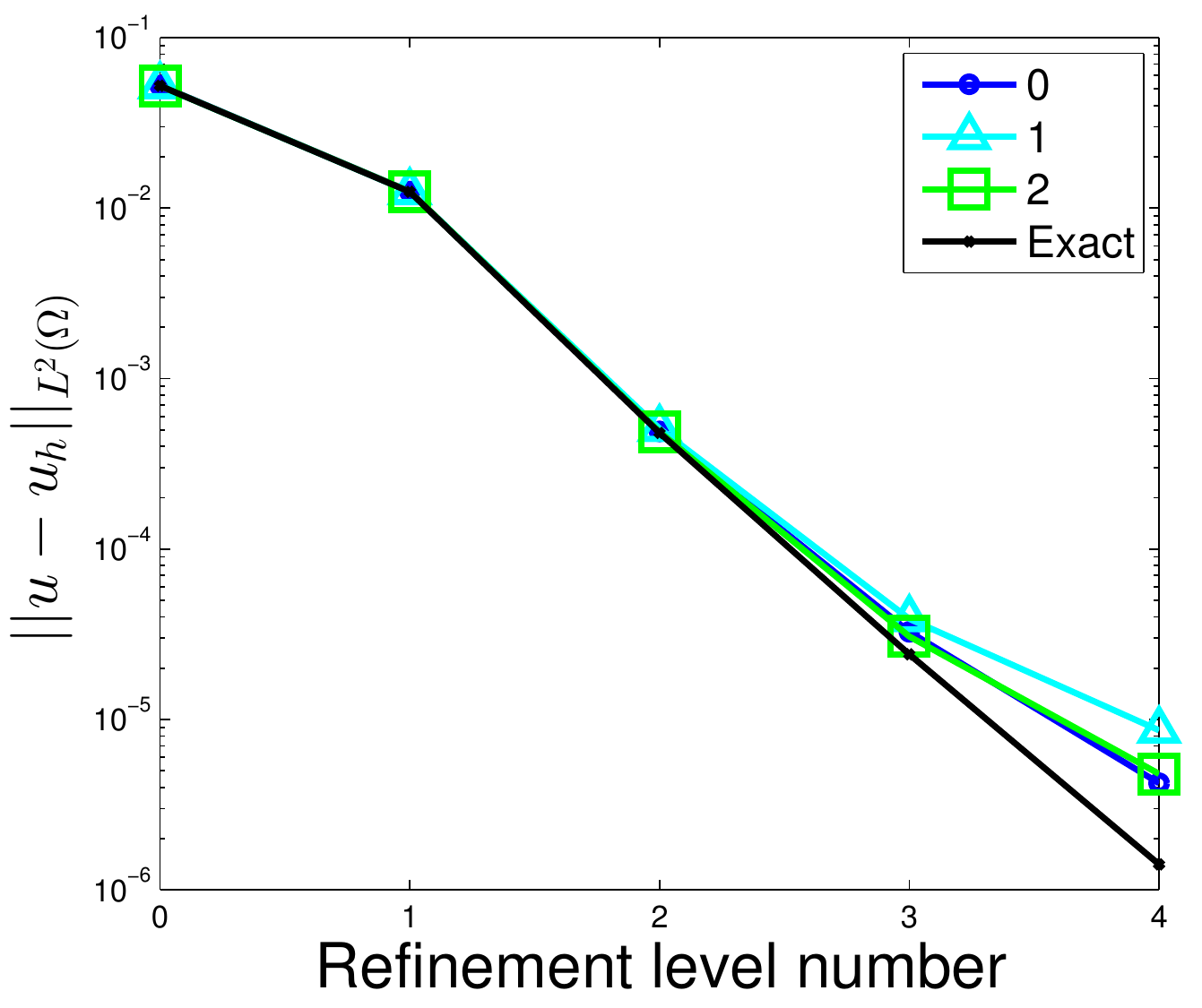}\,\,
\includegraphics[width=.42\textwidth]{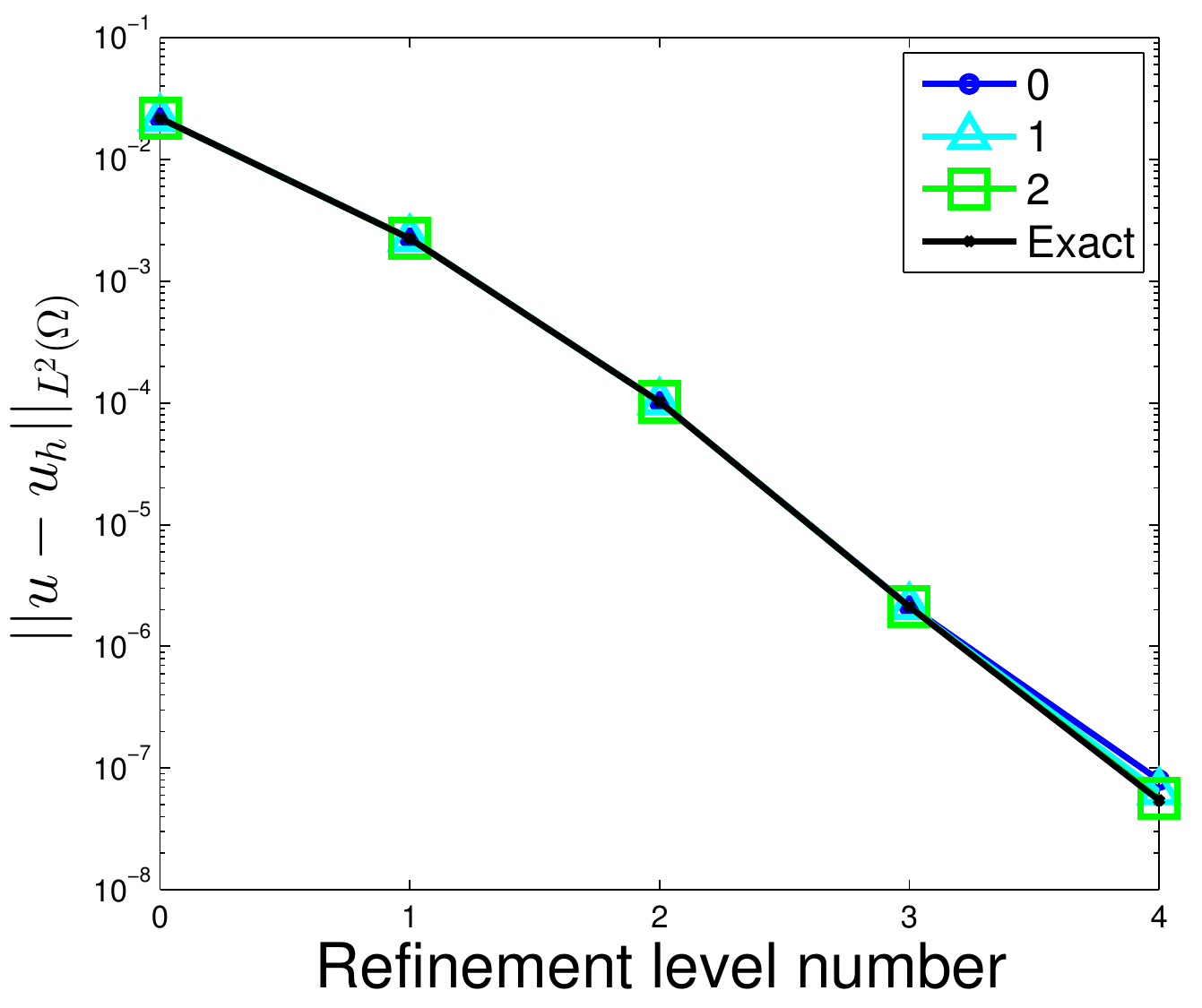}
\end{center}
\caption{3D results - $L^2$ primal error curves for the pairings $P3-P1$ (left)  and $P4-P2$ (right), for the non-symmetric approach and different quadrature rule orders.}\label{fig:3d_pg_31}
\end{figure}

The obtained results are in accordance with the two-dimensional results for both approaches. In Fig.~\ref{fig:3d_slave_p4}, the disturbance for the slave integration approach is shown for the $P4-P4$ pairing. Although not shown here, we note that the results for the $P2-P2$ and $P3-P3$ pairing have a similar behavior.
The non-symmetric approach does not lead to reduced rates considering equal order pairings, i.e., $M_h = M_h^0$, on the refinement levels we considered. As previously, using a lower order dual space, a difference to the exact integration case can be seen. See Fig.~\ref{fig:3d_pg_31} for the disturbance in the primal variable of the $P3-P1$ and $P4-P2$ pairings. 

%% file: figures/axe_picture.tex
\begin{tikzpicture}[scale=0.3]

\tikzstyle{style1} = [black]
\tikzstyle{style2} = [red, very thick, dashed]
\tikzstyle{style3} = [blue, very thick, dash pattern=on \pgflinewidth off 1pt]

\def\a{3.2};
\def\b{1.6};
\def\c{1.6};
\def\d{2.6};

\coordinate (0) at (0, 0);
\coordinate (1) at (\d, -\c);
\coordinate (2) at (0, \a);
\coordinate (3) at (\d, \c);


\draw[thick,->]  (0) -- (1);
\draw[thick,->]  (0) -- (2);
\draw[thick,->]  (0) -- (3);
\def\lnb{0.4};
\def\lna{1};
\def\lnc{1.1};
\node[style1] at ($(\lna,0-\lnc)$) {${\bf x}$};
\node[style1] at ($(\lna,\lnc)$) {${\bf y}$};
\node[style1] at ($(0-\lnb,\b)$) {${\bf z}$};

\end{tikzpicture}

%% file: 05_conclusion.tex
In this article, after reviewing optimal isogeometric mortar methods, a study on the possibility to approximate the mortar integrals by efficient numerical quadrature rules was performed.

To exactly integrate the product of functions defined on non-matching meshes, as in the mortar integrals, it is necessary to construct a merged mesh. Since this construction is of a high complexity, it would be desirable to use a quadrature rule based on the slave mesh only. However, numerical examples show a significant disturbance to the mortar method, especially for higher order splines. Especially the convergence rate of the Lagrange multiplier is reduced to $1/2$ and less.
While  the method improves by increasing the number of quadrature points, the amount of points necessary to obtain nearly optimal results is not predictable. 

To overcome these difficulties, we have considered a non-symmetric saddle point problem based on both master and slave integration rules, which was previously introduced in the finite element context. Numerical examples demonstrate the possibility to reach the accuracy given by an exact integration strategy, although this it is not ensured for all cases.